\documentclass{amsart}

\usepackage{amsmath, amssymb}
\usepackage{mathtools}
\usepackage{enumerate}
\input xy
\xyoption{all}
\numberwithin{equation}{section}
\usepackage{xcolor}

\makeatletter
\@namedef{subjclassname@2020}{%
  \textup{2020} Mathematics Subject Classification}
\makeatother

\begin{document}
\theoremstyle{plain}
\newtheorem{theorem}{Theorem}[section]
\newtheorem{lemma}[theorem]{Lemma}
\newtheorem{corollary}[theorem]{Corollary}
\newtheorem{fact}[theorem]{Fact}

\newtheorem*{theoremA}{Theorem A}
\newtheorem*{theoremB}{Theorem B}
\newtheorem*{theoremC}{Theorem C}
\newtheorem*{theoremD}{Theorem D}

\newtheorem{proposition}[theorem]{Proposition}
\newtheorem{claim}[theorem]{Claim}
\newtheorem*{claim*}{Claim}
\theoremstyle{definition}
\newtheorem{example}[theorem]{Example}
\newtheorem{remark}[theorem]{Remark}
\newtheorem*{remark*}{Remark}
\newtheorem{definition}[theorem]{Definition}
\newtheorem{question}[theorem]{Question}

\def\acvf{\operatorname{ACVF}}
\def\scvf{\operatorname{SCVF}} 
\def\Th{\operatorname{Th}}
\def\tp{\operatorname{tp}}
\def\Int{\operatorname{Int}}
\def\Stab{\operatorname{Stab}}
\def\St{\operatorname{St}}
\def\acl{\operatorname{acl}}
\def\dcl{\operatorname{dcl}}
\def\val{\operatorname{val}}
\def\eq{\operatorname{eq}}
\def\alg{\operatorname{alg}}
\def\aut{\operatorname{Aut}}
\def\trdeg{\operatorname{tr}}
\def\Lin{\operatorname{Lin}}

\renewcommand{\L}{{\mathcal L}}
\renewcommand{\O}{\mathcal O}
\newcommand{\m}{\mathfrak m}

\newcommand{\monster}{\mathcal U}

\def\QQ{\mathbb{Q}}
\def\ZZ{\mathbb{Z}}
\def\CC{\mathbb{C}}
\def\NN{\mathbb{N}}
\def\RR{\mathbb{R}}

\def\an{\mathrm{an}}

\def\bF{\mathbf{F}}
\def\bI{\mathbf{I}}
\def\bZ{\mathbf{Z}}

\newcounter{saveenum}

\newcommand{\elex}{\preccurlyeq}
\newcommand{\elres}{\succcurlyeq}

\newcommand{\todo}[1]{{\bf\small  TODO: #1}}

\def\Ind#1#2{#1\setbox0=\hbox{$#1x$}\kern\wd0\hbox to 0pt{\hss$#1\mid$\hss}
\lower.9\ht0\hbox to 0pt{\hss$#1\smile$\hss}\kern\wd0}
\def\ind{\mathop{\mathpalette\Ind{}}}
\def\notind#1#2{#1\setbox0=\hbox{$#1x$}\kern\wd0\hbox to 0pt{\mathchardef
\nn=12854\hss$#1\nn$\kern1.4\wd0\hss}\hbox to
0pt{\hss$#1\mid$\hss}\lower.9\ht0 \hbox to
0pt{\hss$#1\smile$\hss}\kern\wd0}
\def\nind{\mathop{\mathpalette\notind{}}}

\date{\today}

\title{Definable Equivariant Retractions in Non-Archimedean Geometry}

\author{Martin Hils}
\thanks{MH was partially supported by the German Research Foundation (DFG) via CRC 878, HI 2004/1-1 (part of the French-German ANR-DFG project GeoMod) and under Germany's Excellence Strategy EXC 2044-390685587, `Mathematics M\"unster: Dynamics-Geometry-Structure'.
\\
\indent PS was partially supported by NSF grants no. 1665491 and 1848562.
\\
\indent
All three authors were supported by the Hebrew University of Jerusalem and by the European Research Council under the European Union's Seventh Framework
Programme (FP7/2007--2013)/ERC Grant agreement no. 291111/MODAG}
\address{Institut f\"{u}r Mathematische Logik und Grundlagenforschung, Westf\"{a}lische Wilhelms-Universit\"{a}t M\"{u}nster, Einsteinstr. 62, D-48149 M\"{u}nster, Germany}
\email{hils@uni-muenster.de}

\author{Ehud Hrushovski}
\address{Mathematical Institute, University of Oxford, Andrew Wiles Building, Oxford OX2 6GG, UK}
\email{Ehud.Hrushovski@maths.ox.ac.uk}
\thanks{}

\author{Pierre Simon}
\address{Department of Mathematics, University of California, Berkeley, 733 Evans Hall,
Berkeley, CA 94720-3840, USA}
\email{simon@math.berkeley.edu}
\thanks{}


\begin{abstract}
For $G$ an algebraic group definable over a model of $\acvf$, or more generally a definable subgroup of an algebraic group, we study the stable completion $\widehat{G}$ of $G$, as introduced by Loeser and 
the second author. For $G$ connected and stably dominated, assuming $G$ commutative or that the valued field is of equicharacteristic 0, we construct a pro-definable $G$-equivariant strong deformation retraction of $\widehat{G}$ onto the generic type of $G$.

For $G=S$ a semiabelian variety, we construct a pro-definable $S$-equivariant strong deformation retraction of $\widehat{S}$ onto a definable 
group which is internal to the value group. We show that, in case $S$ is defined over a complete valued field $K$ with value group a subgroup of $\mathbb{R}$, this map descends to an $S(K)$-equivariant 
strong deformation retraction of the Berkovich analytification $S^{\an}$ of $S$ onto a piecewise linear group, namely onto the skeleton of $S^{\an}$. This yields a construction of such a retraction without resorting 
to an analytic (non-algebraic) uniformization of $S$.

Furthermore, we prove a general result on abelian groups definable in an NIP theory: any such group $G$ is a directed union of $\infty$-definable subgroups which all stabilize a generically stable 
Keisler measure on $G$.
\end{abstract}

\maketitle
\setcounter{tocdepth}{1}
\tableofcontents

\section{Introduction}
In \cite{HrLo15}, Loeser and the second author have developed a novel model-theoretic approach to non-Archimedean (algebraic) geometry. For an algebraic variety $V$ defined over some valued field $K$, they construct the so-called \emph{stable completion} of $V$, a pro-definable space, more precisely a functor $\widehat{V}$ which is pro-definable over $K$ and which associates to any valued field extension $L$ of $K$ a 
topological space $\widehat{V}(L)$ naturally containing $V(L)$, endowed with the valuation topology, as a subspace. In case $L$ is non-trivially valued and algebraically closed, $V(L)$ is dense in $\widehat{V}(L)$. 

The construction of $\widehat{V}$ is similar to that of the Berkovich analytification $V^{\an}$ of $V$  (\cite{Ber90}, see also \cite{Duc07}), which is only defined when $K$ is complete with value group  a subgroup of $\RR$. As is the case for $V^{\an}$, the stable completion $\widehat{V}$ has good topological properties, e.g., analogues of local (pathwise) connectedness and local compactness hold for $\widehat{V}$ in the definable category. Assuming that $V$ is quasi-projective, strong topological tameness properties are established for $\widehat{V}$ in \cite{HrLo15}. Most notably, it is shown that $\widehat{V}$ admits a pro-definable strong deformation retraction onto a $\Gamma$-internal space, i.e., onto a piecewise linear space in the definable category. Here, $\Gamma$ denotes the value group. This parallels topological tameness results in the case of $V^{\an}$, established by Berkovich. Actually, assuming $K$ is complete with $\Gamma_K\leq\mathbb{R}$, it is shown in \cite[Chapter 14]{HrLo15} that, for some suitable $K^{max}\supseteq K$, $V^{\an}$ may be obtained as a topological quotient of $\widehat{V}(K^{max})$ and the pro-definable strong deformation retraction descends to a strong deformation retraction of $V^{\an}$ onto a finite simplicial complex, thus reproving and generalizing results by Berkovich, as this holds for every quasi-projective variety $V$, without any smoothness assumption on $V$. Moreover, in this context, local contractibility of $V^{\an}$ is established in \cite{HrLo15}, as well as the finiteness of the 
number of homotopy types in $\{\mathcal{V}_b^{\an}\mid b\in S(K)\}$, where $\mathcal{V}\rightarrow S$ is a quasi-projective family of algebraic varieties defined over $K$. For a survey of the results in \cite{HrLo15} and in particular the consequences in the realm of Berkovich spaces, we refer to \cite{Duc13}.

The main tools used in \cite{HrLo15} concern the geometric model theory of the (first order) theory $\acvf$ of algebraically closed  non-trivially valued fields. The study of $\acvf$ goes back to Abraham Robinson who established its model-completeness. His arguments actually yield quantifier elimination for $\acvf$ in various natural languages. Haskell, Macpherson and the second author initiated the study of $\acvf$ from a point of view of geometric model theory. In \cite{HHM}, they show that the imaginary sorts in $\acvf$ are classified by certain natural group quotients, the so-called \emph{geometric sorts}, and they identify the stable part of $\acvf$ as the collection of definable sets which are internal in the residue field $k$. 
In \cite{HHM2} the authors systematically develop stable domination as a means to lift phenomena known in stable theories to the unstable context. Stably dominated types are definable types whose generic extension is entirely controlled by its stable part. In $\acvf$, these correspond to the definable types orthogonal to $\Gamma$.

Let us briefly describe the construction of the Berkovich analytification of an algebraic variety. Suppose that $K$ is a field which is complete with respect to a non-Archimedean norm $\mid\!\cdot\!\mid:K\rightarrow\RR_{\geq0}$ and that $V$ is an 
algebraic variety defined over $K$. The construction of the anaytification $V^{\an}$ is done by gluing affine pieces, so let us assume that $V$ is affine. Then $V^{\an}$ is  the set of all multiplicative semi-norms $\mid\!\cdot\!\mid_v:K[V]\rightarrow\RR_{\geq0}$ extending $\mid\!\cdot\!\mid$ on $K$, and the topology on $V^{\an}$ is defined as the weakest topology such that for any $f\in K[V]$, the map $\mid\!\cdot\!\mid_v\mapsto \mid\! f\!\mid_v$ is continuous. In this way, the natural embedding 
of $V(K)$ into $V^{\an}$ becomes a homeomorphism. 

Identifying $\RR_{\geq0}$ with $\RR\cup\{\infty\}$ via $-\log$, the field $K$ becomes a (complete) valued field with $\Gamma_K\leq\RR$, and, by quantifier elimination in $\acvf$, as a set, $V^{\an}$ equals  
$$\{p=\tp(\overline{a}/K\cup\RR)\in S_V(K\cup\RR)\mid\Gamma_{K(\overline{a})}\leq\RR\},$$
where $S_V(K\cup\RR)$ denotes the set of complete types concentrating on $V$, in the theory $\acvf$, over the (2-sorted) parameter set $K\cup\RR$. 

The key insight of \cite{HrLo15} is that one may obtain a model-theoretic analogue of $V^{\an}$ 
by considering those (global definable) types concentrating on $V$ whose realizations do not increase the value group, relative to the base model over which one works. Consequently, the stable completion $\widehat{V}$ is defined as the set of stably dominated types 
concentrating on $V$, endowed with a topology which is similarly defined as in the case of $V^{\an}$. It is an important feature that, when one identifies a definable type with its canonical base, $\widehat{V}$ naturally gets the structure of a (strict) pro-$K$-definable space in $\acvf$. Thus model-theoretic methods may be applied to $\widehat{V}$, 
in particular powerful tools both from stability theory  and from $o$-minimality, two well developed and completely different strands of geometric model theory. Any type $\tp(\overline{a}/K)$ with $\trdeg(K(\overline{a})/K)=\trdeg(k_{K(\overline{a})}/k_K)$ is stably dominated, and the collection of these types, which corresponds to the sub-collection of all strongly stably dominated types, forms an ind-definable subspace $V^{\#}$ of $\widehat{V}$.

Now let us move to the equivariant situation. If $K$ is a complete field with respect to an non-Archimedean norm and $G$ is an algebraic group defined over $K$, then $G(K)$ naturally acts on $G^{\an}$, and one may ask whether there is a $G(K)$-equivariant 
retraction of $G^{\an}$ onto a skeleton which is then naturally a piecewise linear group. For $E$ a Tate elliptic curve defined over $K$, i.e., an elliptic curve with bad (split) multiplicative reduction, $E^{\an}$ does indeed admit an $E(K)$-equivariant strong deformation retraction onto the circle group $\RR/\ZZ$, as one may for example show using the Tate uniformization; in case $E$ is an elliptic curve with good reduction, $E^{\an}$ equivariantly retracts to the trivial group. Actually, whenever $S$ is an abelian or even semi-abelian variety defined over $K$, passing to a finite separable extension of $K$ if necessary, $S^{\an}$  admits an analytic uniformization by an analytic group which is equivariantly contractible, from which one obtains an 
$S(K)$-equivariant strong deformation retraction of $S^{\an}$ onto its skeleton  (see \cite[Section~6.5]{Ber90}, where this is explained for abelian varieties).

In the present paper, we prove analogous results for the stable completion $\widehat{S}$ of a semiabelian variety $S$ defined over some valued field $K$, confirming that the model-theoretic approach to non-Archimedean geometry developed in 
\cite{HrLo15} is very adequate for the study of semiabelian varieties as well, as the topological tameness properties one expects do also hold in this context. We would like to highlight that, contrarily to the classical approach via analytic uniformizations (see \cite{BoLu84}), our construction may be performed over an arbitrary valued field (with value group not necessarily Archimedean) and does not require the passage to a finite separable extension in the process. 
Moreover, we reprove the results for $S^{\an}$ mentioned in the previous paragraph 
in the same way topological tameness results for $V^{\an}$ could be deduced in \cite{HrLo15} from the corresponding results for the stable completion $\widehat{V}$ in the definable category.

Before we may state the main contributions of our paper, we need to mention a result from \cite{HrRi19}, which is a crucial ingredient in our construction.

In \cite{HHM2}, an important and very useful  structural way to decompose types in $\acvf$ is obtained: if $M$ is a maximally complete model of $\acvf$ and $\overline{a}$ is a (possibly imaginary) tuple from the monster model $\monster$, 
then $\tp(\overline{a}/M\cup\Gamma(M\overline{a}))$ is stably dominated. So types may be understood, to some extent, by types in the value group and types in the residue field.\footnote{This may be seen as a very powerful incarnation of the 
Ax-Kochen-Ershov philosophy which postulates that model-theoretic questions about henselian valued fields may be reduced to questions about the residue field and the value group.} Moreover, $\acvf$ has the invariant extension property, i.e., any type over an algebraically closed set admits a global automorphism 
invariant extension. It follows that $\acvf$ is a \emph{metastable} theory. 

Metastability is an axiomatic framework capturing the phenomena we described in the previous paragraph. It  has been introduced in \cite{HHM2} and further developed in work of Rideau-Kikuchi and the second author \cite{HrRi19},  who undertake a thorough study of groups definable in metastable theories. In their work, stably dominated groups, i.e., ($\infty$-)definable groups admitting a stably dominated generic type, play a key role, an important reason being that many features of stable groups lift to them.
As one of the main results, they establish a group version of the metastability property in the commutative case: for any definable abelian group $G$ in $\acvf$, there is a definable homomorphism $\lambda:G\twoheadrightarrow\Lambda$ such that $\Lambda$ is $\Gamma$-internal and $N:=\ker(\lambda)$ is a connected group which is 
an increasing union (indexed essentially by the value group) of definable stably dominated subgroups.

We show that whenever $G=S$ is a semiabelian variety, $N$ is in fact itself stably dominated. More precisely, if $S$ is defined over a valued field $F\subseteq K\models\acvf$, there is an $F$-definable decomposition \begin{equation}0\rightarrow N\rightarrow S\rightarrow\Lambda\rightarrow0\label{eq:Decomp}\end{equation} such that $\Lambda$ is $\Gamma$-internal and $N$ stably dominated, definable and connected. In particular, 
it follows that $p_N\in\widehat{S}$, where $p_N$ denotes the generic type of $N$.

For abelian varieties, stable domination of $N$ is shown in \cite{HrRi19}. It is easy to see that this also holds for an algebraic torus, e.g., if $G=\mathbb{G}_m^n$, then $N=(\mathcal{O}^*)^n$ is stably dominated, with $G/N\cong\Gamma^n$. We show that stable domination of the corresponding $N$ is preserved in short exact sequences 
of algebraic groups. As a semiabelian 
variety is an extension of an abelian variety by an algebraic torus, we may thus lift the result to the semiabelian case. 

For $G$ any definable group, we denote by $p_e\in\widehat{G}$ the type of the identity element in $G$. We obtain the following result.

\begin{theoremA}[Theorem~\ref{T:Retract-all-char}]
Let $S$ be a semi-abelian variety defined over $F\subseteq K\models\acvf$, and let $0\rightarrow N\rightarrow S\rightarrow\Lambda\rightarrow0$ be the decomposition from \eqref{eq:Decomp}. 

Then there is an $F$-definable special\footnote{See Definition~\ref{D:SpecialRetr} for the notion of a special deformation retraction.} deformation retraction $\rho: [0,\infty] \times \widehat S \to \widehat S$ with final image $\Sigma \subseteq \widehat{S}$ such that $\rho$ is equivariant under the action of $S$ by multiplication and  for each $t<\infty$, $q_t=\rho(t,p_e)$ is the generic type of a connected strongly stably dominated definable subgroup $N_t$ of $N$, with $N_t$ Zariski dense in $S$.

Moreover, the morphism $\pi:S\rightarrow\Lambda$ induces definable bijection between $\Sigma$ and $\Lambda$.
\end{theoremA}

Modulo some continuity issues concerning the map given by the tensor product $\otimes$ (dealt with in Subsection~\ref{Sub:ContTensor}), Theorem~A follows in a rather straight forward way from the existence of a definable continuous path in $\widehat{N}$ between $p_e$ and $p_N$ (the generic type of $N$), along generic types of (Zariski-dense strongly stably dominated) definable subgroups $N_t$ of $N$. Let us illustrate this first with two easy examples. Denote by $\eta_{c,\gamma}$ the generic type of the closed ball $B_{\geq\gamma}(c)$, for $\gamma\in\Gamma\cup\{\infty\}$. 

\begin{itemize}
\item The standard deformation retraction $\rho:[0,\infty]\times\widehat{\mathcal{O}}\rightarrow\widehat{\mathcal{O}}$, sending $(\gamma,\eta_{c,\delta})$ to  $\eta_{c,\min(\delta,\gamma)}$  is definable and $(\mathcal{O},+)$-equivariant with final image  $\{\eta_{0,0}\}$, and as an equivariant continuous map it is entirely determined by the path $$q:[0,\infty]\rightarrow\widehat{\mathcal{O}}, \,\,q(\gamma):=\rho(\gamma,p_e)=\eta_{0,\gamma},$$ where $p_e=\eta_{0,\infty}$. Thus, $N_\gamma=\gamma\O$ in this example.
\item The map $\rho':[0,\infty]\times\mathbb{G}_m\rightarrow\widehat{\mathbb{G}_m}$, $(\gamma,c)\mapsto\eta_{c,v(c)+\gamma}$ extends uniquely to a $\mathbb{G}_m$-equivariant deformation retraction $\rho:[0,\infty]\times\widehat{\mathbb{G}_m}\rightarrow\widehat{\mathbb{G}_m}$, via
$$\rho(\gamma,\eta_{c,v(c)+\delta})=\eta_{c,v(c)+\min(\gamma,\delta)} \text{\,\,\,(for $c\neq0$ and $\delta\geq0$)}.$$
Its final image is $\{\eta_{c,v(c)}\!\mid \!c\neq0\}=\{\eta_{0,\gamma}\mid\gamma\in\Gamma\}\cong\Gamma$, and one has $p_e=\eta_{1,\infty}$, $N=N_0=\O^*$ and $N_\gamma=1+\gamma\O$ for all $\gamma\in (0,\infty)$. Setting $q_\gamma=\rho(\gamma,p_e)=\eta_{1,\gamma}$,  one may check that $$\rho(\gamma,p)=\widehat{m}(q_\gamma\otimes p),$$
the convolution of $q_\gamma$ and $p$. Here, $m$ denotes the multiplication in $\mathbb{G}_m$.
\end{itemize}

\smallskip

In our paper, we give two proofs of the existence of a group path as above, one in arbitrary characteristic which is valid in all stably dominated abelian 
definable subgroups of an algebraic group, one in equicharacteristic 0, valid for all stably dominated subgroups of an algebraic group:

\begin{theoremB}[Theorem~\ref{P:Group-path-general} \& Theorem~\ref{P:continuity}]
Let $G$ be an algebraic group defined over a valued field $F$, and let $N$ be an $F$-definable stably dominated connected subgroup of $G$. Assume that $N$ is commutative or that $F$ is of equicharacteristic 0. Then there is an $F$-definable path 
$q:[0,\infty]\rightarrow N^{\#}$ such that 
\begin{enumerate}
\item [(i)] $q_{\infty}=p_e$ and $q_{0}=p_N$; 
\item[(ii)] for any $t<\infty$, $q_t$ is the generic type of a connected stably dominated definable subgroup $N_t$ of $N$ which is Zariski dense in $N$.
\end{enumerate}
\end{theoremB}

The construction in arbitrary characteristic relies on the existence of a definable path between $p_e$ and $p_N$, which is then turned into a group path, using an averaging process (described in \cite{HrRi19} and requiring commutativity of the group, for technical reasons) reminiscent of Zilber's Indecomposability Theorem. The second construction, only valid in equicharacteristic 0, is more explicit and does not require the group $N$ to be commutative. It utilizes an intrinsic scale provided by the value group, and the subgroup $N_t$ in the group path is found 
as the kernel of the homomorphism to the maximal $\O/t\O$-internal quotient of $N$. The proof may be understood as a linearization procedure at the level of generic types, which is available even if, in $\acvf$, there exists no exponential map in the definable category.

\smallskip

As we have already mentioned, we deduce Theorem~A from Theorem~B. Moreover, Theorem~B yields the following strong contractibility result as a corollary.

\begin{theoremC}[Theorem~\ref{T:Contraction-gp} \& Corollary~\ref{C:Contract-char-00}]
Let $N$ be a group which satisfies the assumptions of Theorem~B. Then $\widehat{N}$ admits an $F$-definable $N$-equivariant special deformation retraction with final image $\{p_N\}$.
\end{theoremC}

\smallskip
It is worth mentioning the following connection to group schemes over $\mathcal{O}$. Assume the base field $F$ is algebraically closed. 
If $H$ is a connected Zariski-dense stably dominated ($\infty$-definable) subgroup of the affine algebraic group $G$, there is a group scheme $\mathbb{H}$ over $\mathcal{O}$ such that $H$ is 
definably isomorphic to the group of $\mathcal{O}$-valued points of $\mathbb{H}$. (See \cite[Theorem~6.11]{HrRi19}.) Under additional assumptions on $\mathbb{H}$, a chain of stably dominated subgroups becomes visible geometrically, namely via the map $\val (c) \mapsto \ker(\mathbb{H}(\mathcal{O}) \to \mathbb{H}(\mathcal{O}/c\mathcal{O}))$.

\smallskip

Assume now that $F$ is a valued field with value group $\Gamma_F\leq\RR$ and that $V$ is an algebraic variety defined over $F$. In \cite[Chapter~14]{HrLo15}, strong topological links are established between the stable completion $\widehat{V}$ of $V$ (the \emph{definable} world) and the Berkovich analytification $V^{\an}$ (the \emph{analytic} world). If $F^{\max}$ 
denotes an algebraically closed valued field extension of $F$ with value group $\Gamma_{F^{\max}}=\RR$ and residue field $k_{F^{\max}}=(k_F)^{\alg}$, then there is a natural continuous surjective  and closed map $\pi_V:\widehat{V}(F^{\max})\rightarrow V^{\an}$. Moreover, considering the parameter set  $\bF:=(F,\RR)$, any pro-$\bF$-definable (special) 
deformation retraction $$H:[0,\infty]\times\widehat{V}\rightarrow\widehat{V}$$ descends to a (special) deformation retraction $$\widetilde{H}:\RR^+_{\infty}\times V^{\an}\rightarrow V^{\an},$$ where $\RR^+_\infty=[0,\infty](\bF)=\RR_{\geq0}\cup\{\infty\}$. If the final image $Z$ of $H$ is $\Gamma$-internal,  then the final image of $\widetilde{H}$ equals 
${\bf{Z}}=\pi_V(Z(F^{\max}))$, a set which is naturally homeomorphic to a piecewise linear subset of $\RR_{\infty}^n$ (where $\RR_\infty=\RR\cup\{\infty\}$ and $n\in\NN$), carrying a $\QQ$-tropical structure\footnote{This means that all the coefficients in the inequalities describing $\bf{Z}$ may be taken from $\QQ$.}.

We call such deformations $\widetilde{H}$ \emph{definably induced}. In case $V=G$ is an algebraic group and $H$ is $G$-equivariant, the map $\widetilde{H}$ will automatically be $G(F)$-equivariant. Theorem~A thus yields the following topological tameness result about analytifications of semi-abelian varieties.

\begin{theoremD}[Theorem~\ref{T:Berk-contract}]\label{T:D}
Let $S$ be a semi-abelian variety defined over a valued $F$ with $\Gamma_F\leq\RR$. Consider the $F$-definable decomposition $0\rightarrow N\rightarrow S\rightarrow\Lambda\rightarrow0$ from \eqref{eq:Decomp}. 
Then there is a definably induced  $S(F)$-equivariant special deformation retraction $$\tilde{\rho}: \RR^+_\infty\times S^{\an}\rightarrow S^{\an}$$ with final image a skeleton $\bf{\Sigma}$. The map $\pi:S\rightarrow\Lambda$ induces a bijection between $\Sigma$ and  $\Lambda(\bF)$, where $\Lambda(\bF)$ is a subset of $\RR^n$ carrying the structure of a piecewise linear abelian group.
\end{theoremD}

\smallskip
Let us mention that in case $A$ is an abelian variety over a valued field $F$ with $\Gamma_F\cong\mathbb{Z}$ and if $N\leq A$ is as in \eqref{eq:Decomp}, then $N(F)$ corresponds to the $\mathcal{O}_F$-points of the identity component of the N\'eron model $\mathcal{A}$ of $A$.\footnote{We thank Antoine Ducros and Peter Schneider for pointing this out.} Observe, though, that the (maximal) stably dominated subgroup $N$ of $A$ always exists, even when $A$ is definable over a valued field with non-Archimedean value group.

\smallskip

In the present article, even if we recall the relevant key concepts and results, we assume the reader is familiar with the model theory of $\acvf$ as developed in \cite{HHM,HHM2} and more specifically with the model-theoretic approach to non-Archimedean geometry 
from \cite{HrLo15}. Moreover, we will use  results from \cite{HrRi19} on stably dominated groups at many places, and so some acquaintance with these is certainly helpful.

\smallskip

Here is an overview of the article. In Section~\ref{S:Prelim}, we gather a number of facts which we will use 
in our paper, in particular on stable domination, the stable completion of an algebraic variety and on stably dominated groups in $\acvf$. 
The material we present is mostly from \cite{HHM2,HrLo15,HrRi19}, but some facts are new, in particular one on the maximal internal quotient of a generically stable group (Proposition~\ref{P:metastable-group}) and a continuity 
result concerning $\otimes$ in the stable completion (Proposition~\ref{P:continuity-tensor}), which are both key ingredients in our construction.

In Section~\ref{S:AllChar}, we give our first construction of a definable group path, valid in all characteristics, and in Section~\ref{S:Retract0}, we present the explicit construction in equicharacteristic 0. The results about equivariant retractions are 
then transferred to the setting of Berkovich analytifications in Section~\ref{S:Berkovich}, where we prove in particular Theorem~D. Finally, in Section~\ref{S:NIP-abgp}, we prove a general result on commutative groups definable in an arbitrary NIP theory (Proposition~\ref{P:stab-Keisler}), which is a rather weak analogue of our main result, namely that any such group $G$ may be written as a directed union of $\infty$-definable subgroups which all stabilize a generically stable 
Keisler measure on $G$.

\section{Tying up some loose ends}\label{S:Prelim}
We work in a complete theory $T$, and $\monster$ denotes a monster model of $T$. We assume throughout that $T$ eliminates quantifiers and imaginaries. Our notation and terminology is mostly standard. If $D$ is a $C$-definable set, for $C$ a subset of $\monster$, we denote by $D(C)$ the set $D(\monster)\cap\dcl(C)$. We will use various generalizations of definable sets in our paper, in particular (strict) pro-definable and $\infty$-definable sets, relatively definable subsets of such sets, and occasionally ind-definable and iso-definable sets. See \cite[Section 2.2.]{HrLo15} or \cite{HrRi19} for definitions and properties of these notions. For basic facts about the model theory of $\acvf$, we refer to \cite{HHM,HHM2}.

\subsection{Stably embedded sets}
Let $C$ be a small subset of $\monster$. Recall that a (small) family $(X_i)_{i\in I}$ of $C$-definable sets in $\monster$ is \emph{stably embedded} if for any finite sequence $(i_1,\ldots,i_m)$ of elements of $I$, any 
$\monster$-definable subset of $\prod_{j=1}^mX_{i_j}$ is definable with parameters from $C\cup\bigcup_{i\in I}X_i$. If $X$ and $D$ are definable sets, $X$ is said to be \emph{$D$-internal} if ($X=\emptyset$ or) there 
are $m\in\NN$ and a surjective $\monster$-definable function $g:D^m\rightarrow X$.

\emph{Generically stable} types are definable types that behave in some ways like types in a stable theory. For a definition of generic stability in an arbitrary theory $T$, see, e.g., \cite{AdCaPi14}. We will only be concerned with generically stable types in NIP theories, where there are various useful characterizations of this notion (see \cite[Theorem~2.29]{NIPbook}).
\begin{lemma}\label{L:intgerm}
Let $p$ be a $C$-definable generically stable type, and $D$ a $C$-definable set. Let $f_{-}$ be a definable family of functions to $D$. Then the set of $p$-germs of instances of $f_{-}$ is $D$-internal.
\end{lemma}

\begin{proof}
Let $X$ be the set of $p$-germs of instances of $f_{-}$, so $X$ is a definable set since the type $p$ is definable. Let $a_1,\ldots,a_n$ be a sufficiently long Morley sequence of $p$ over $C$. For a parameter $c$, let $g(c)=(f_c(a_1),\ldots,f_c(a_n))\in D^n$. Taking $n$ large enough, if $g(c)=g(c')$  the germs of $f_c$ and $f_{c'}$ are equal. Hence the map $\pi$ from the image of $g$ to $X$ sending $g(c)$ to the $p$-germ of $f_c$ shows that $X$ is $D$-internal.
\end{proof}

\begin{lemma}\label{L:intstemb}
($T$ is any theory.) Let $X$ be an $\emptyset$-definable set. Assume that $X$ is stably embedded, then so is $\Int(X)$: the union of the $\emptyset$-definable $X$-internal sets.
\end{lemma}
\begin{proof}
First, let $Y$ be $\emptyset$-definable and $f_a: X \to Y$ a definable bijection. Consider the $\emptyset$-definable set $A$ of codes of bijections $f_{a'}:X\to Y$. For $f\in A$, define $g = f^{-1}\circ f_{a}$. Then $g : X\to X$ is a bijection. As $X$ is stably embedded, $g = h_b$ for some definable function $h_b$, with $b$ in some $X^n$. The function mapping $b$ to $f$ such that $h_b =  f^{-1}\circ f_{a}$ is $a$-definable and defines a surjection from some definable subset of $X^n$ to $A$. Hence $A$ is $X$-internal. Now any definable subset of $Y$ is definable using a parameter from $A$ and parameters from $X$ (since $X$ is stably embedded), therefore any definable subset of $Y$ is definable with parameters from $\Int(X)$.

Now if $Y$ is an $\emptyset$-definable $X$-internal set, there is a definable surjection $f_a : X^m \to Y$. Let $A$ be the set of codes of surjections $f_{a'}: X^m \to Y$. Then $A$ is $\emptyset$-definable and any definable subset of $Y$ is definable with parameters in $X$ and $A$. It remains to show that $A$ is $X$-internal. The argument is the same as above: let $f,f' \in A$, then we can define the correspondence $f'^{-1}\circ f := \{(x,x')\in X^m \times X^m : f(x)=f'(x')\}$. By stable embeddedness of $X$, any such correspondence can be written as $h_b$ for some $b\in X^k$. Now fixing an element $f\in A$, we see that any element of $A$ is definable from such a $b\in X^k$ along with $f$, hence $A$ is $X$-internal as required.
\end{proof}


\subsection{Stably dominated and strongly stably dominated types} We briefly recall some facts around (strongly) stably dominated types. Given a set of parameters $C$, we denote by $\St_C$ the collection of all $C$-definable stable stably embedded sets. For a tuple $a$, we denote by $\St_C(a)$ the set $\dcl(Ca)\cap \St_C(\mathcal{U})$. 

In the following definition, 
non-forking $\ind$ is meant with respect to the stable multi-sorted structure $\St_C(\mathcal{U})$.

\begin{definition}[\mbox{\cite{HHM2,HrRi19}}]
A type $p=\tp(a/C)$ is called \emph{stably dominated} if for any tuple $b$ with $\St_C(a) \ind_{\St_C(C)} \St_C(b)$, one has $\tp(a/\St_C(a))\vdash\tp(a/\St_C(a)b)$  (or equivalently, $\tp(b/\St_C(a))\vdash\tp(b/Ca)$). It is called 
\emph{strongly stably dominated} if there exists a formula $\phi(x)\in \tp(a/\St_C(a))$ such that for any tuple $b$ with $\St_C(a) \ind_{\St_C(C)} \St_C(b)$, one has $\phi(x)\vdash\tp(a/\St_C(a)b)$. 
\end{definition}

Clearly, every strongly stably dominated type is stably dominated.


By \cite[Proposition~2.6.12]{HrLo15}, if $p$ is a definable type based on some $C=\acl(C)$ and $p|A$ is strongly stably dominated for some $A=\acl(A)$, then $p|C$ is strongly stably dominated.

\begin{fact}[\cite{HHM2}] Let $p$ be a global definable type in some theory $T$.
\begin{enumerate}
\item If $p$ is stably dominated, then it is generically stable.
\item In $T=\acvf$, $p$ is stably dominated if and only if it is generically stable if and only if it is orthogonal to $\Gamma$.\footnote{Recall that a global definable type $p$ is called \emph{orthogonal to $\Gamma$} if $\Gamma(\monster a)=\Gamma(\monster)$ for $a\models p$.}
\end{enumerate}
\end{fact}

The following is a special case of \cite[Proposition~8.1.2]{HrLo15}.
\begin{fact}\label{fact:ssd}
We work in $\acvf$. Let $q$ be an $A$-definable type on a variety $V$. Then $q$ is strongly stably dominated if and only if $\dim(q)=\dim(h_*q)$ for some $A$-definable map $h$ into a stable sort; where $\dim(h_*q)$ refers to Morley dimension.
\end{fact}

\subsection{The stable completion of an algebraic variety}We work in $\acvf$, and we will now briefly recall the notion of the \emph{stable completion} of an algebraic variety from  \cite{HrLo15}.\footnote{We refer to \cite{HrLo15} and also to \cite{Duc13} for details.} Given an algebraic variety $V$ defined over a valued field $F$, its stable completion $\widehat{V}$ is a functor which to every set $B$ over which $V$ is defined 
associates the set $\widehat{V}(B)$ of all global $B$-definable stably dominated types concentrating on $V$. It is naturally given by a strict pro-$F$-definable set, i.e., by a projective limit of $F$-definable sets with surjective transition functions.

One endows $\widehat{V}$ with the following ind-definable topology: for any $F\subseteq K\models\acvf$, a basis of open subsets for $\widehat{V}(K)$ is given by sets of the form 
$\{p\in\widehat{U}\mid g_*(p)\in \Omega\}$, for $U$ an open affine subvariety of $V$ defined over $K$, $\Omega\subseteq\Gamma_\infty$ an open subset and $g=\val\circ G$ for some regular function $G$ on $U$. Here, $g_*(p)\in\Gamma_\infty(K)$ is defined as $g(a)$, where $a\models p|K$. As $p$ is orthogonal to $\Gamma$, this is well-defined.

If $X\subseteq V$ is a definable subset, $\widehat{X}:=\{p(x)\in\widehat{V}\mid p(x)\models x\in X\}$ is a relatively definable subset of $\widehat{V}$ which is endowed with the subspace topology. 
One denotes by $X^{\#}(B)$ the set of all strongly stably dominated types in $\widehat{X}(B)$. Then $X^{\#}$ is naturally an ind-definable set, and we endow it with the subspace topology.

We refer to \cite{HrLo15} for the definition and study of definable analogues of various topological properties, in particular definable compactness and definable (path) connectedness.

\subsection{Pro-definable groups and generic types}
We now recall the concept of genericity in the context of pro-definable groups, following \cite[Section~3]{HrRi19}. Let $G$ be a pro-$C$-definable group,\footnote{See \cite{HrRi19} for the definition of pro-definable and $\infty$-definable groups.} where $C$ is some small parameter set. A global $C$-definable type $p$ concentrating on $G$ is called \emph{right generic} in $G$ over 
$C$ if for all $g\in G(\mathcal{U})$, the type $g\cdot p$ is $C$-definable. Left genericity is defined in a similar way, as is right genericity in a pro-definable principal homogeneous space under $G$. A right or left generic type $p$ of $G$ is called \emph{symmetric} if for any other global definable type $q$ concentrating on $G$ one has $p\otimes q=q\otimes p$. Note that in an NIP theory, this is equivalent to generic stability of $p$.

If $G$ admits a smallest pro-$\monster$-definable subgroup $H$ of bounded index, $H$ is called the \emph{strong connected component} of $G$ and denoted by $G^{00}$. Similarly, if the intersection $H$ of all pro-$\monster$-definable subgroups of finite index in $G$ has bounded index in $G$, then $H$ is called the \emph{connected component} of $G$ and denoted by $G^{0}$. 
The group $G$ is called \emph{connected} if $G=G^{0}$, and \emph{strongly connected} if $G=G^{00}$.

The following fact combines \cite[3.4, 3.9, 3.11 \& 3.12]{HrRi19}.

\begin{fact}($T$ is any theory.)  Let $G$ be a pro-definable group.\label{F:Generics}

Assume that $p$ is a right (resp.\ left) generic type of $G$. Then the following holds:
\begin{enumerate}
\item $G$ is pro-definably isomorphic to a pro-limit of definable groups. In particular, if $G$ is  $\infty$-definable, it is an intersection of definable groups. \label{I:ProGroup}
\item $G^{00}$ exists and one has $G^{00}=G^0=\mathrm{Stab}(p)$, where $\mathrm{Stab}(p)$ denotes the left (resp.\ right) stabilizer of $p$.   \label{I:Connected-component}
\setcounter{saveenum}{\value{enumi}}
\end{enumerate}
Assume in addition that $p$ is symmetric. Then the following holds:
\begin{enumerate}
\setcounter{enumi}{\value{saveenum}}
\item Right and left generic types coincide in $G$, they are all symmetric, and $G$ acts transitively on the set of generic types, which is bounded. \label{I:SymGen}

In particular, $p$ is left and right generic, and 
$G$ is connected if and only if $p$ is the unique (right) generic type of $G$. 
\end{enumerate}
\end{fact}

In our work, mostly symmetric (right) generic types will play a role. We will call them \emph{generic} types in what follows.

\begin{definition}A pro-definable group is called \emph{generically stable} (\emph{stably dominated} or \emph{strongly stably dominated}, respectively) if it admits a generic type which is generically stable (stably dominated or strongly stably dominated, respectively).
\end{definition}

Note that all parts of Fact~\ref{F:Generics} apply to all generically stable and thus in particular to all stably dominated groups.

\begin{lemma}\label{lem:coset-generic}
($T$ is any theory.) Let $G$ and $H\leq G$ be pro-definable groups. Assume that $G$ is stably dominated and that $H$ admits a right generic type. Let $a$ be generic in $G$. Then $a$ is right generic in the coset $a\cdot H$ over $\pi(a)$, where $\pi:G\rightarrow G/H$ is the canonical projection.
\end{lemma}

\begin{proof}
By Fact~\ref{F:Generics}, we may assume that $a\in G^{0}$. Let $b$ be right generic in $H$ over $a$, such that $b\in H^{0}$, and so $b\in G^{0}$ as well. By symmetry of the generic type of $G$, $a$ is generic in $G$ over $b$, as is $ab^{-1}$. By symmetry again, $b$ is right generic in $H$ over $ab^{-1}$. Being a translate of $b$ over $ab^{-1}$,  $a = (ab^{-1})\cdot  b$  is right generic in $a\cdot H$ over $ab^{-1}$. (Note that $a\cdot H$ is definable over $ab^{-1}$.)
\end{proof}

The following proposition is a key ingredient for our second construction of definable equivariant retractions in Section~\ref{S:Retract0}.

\begin{proposition}\label{P:metastable-group}
Let $G$ be a pro-$C$-definable generically stable connected group, and let $D$ be a stably embedded $C$-definable set. There exists a pro-$C$-definable group $\mathfrak g_D$ internal to $D$ and a pro-$C$-definable homomorphism $g:G\to \mathfrak g_D$ which is maximal in the sense that for any other pro-$C$-definable group homomorphism $g':G\to \mathfrak g'_D$ with $\mathfrak g'_D$ internal to $D$, $g'$ factors through $g$.

The generic of $\mathfrak g_D$ is interdefinable over $C$ with $\dcl(Ca)\cap \Int_C(D)$, where $a$ is a generic of $G$ over $C$ and $\Int_C(D)$ is the union of $C$-definable, $D$-internal sets.
\end{proposition}

\begin{proof}
We follow the proof of \cite[Proposition~4.6]{HrRi19}  very closely. By Lemma~\ref{L:intstemb}, $\Int_C(D)$ is stably embedded. Let $p$ be the (by Fact~\ref{F:Generics} unique) generic type of $G$. For any $a$, let $\theta(a)$ enumerate $\dcl(Ca)\cap \Int_C(D)$. Fixing $a\in G$ generic over $C$, consider the map $f_a: b \mapsto \theta(ab)$ which is a pro-$C$-definable map on $G$. As the generic $p$ of $G$ is generically stable, the $p$-germ $\tilde a$ of $f_a$ is strong by \cite[Theorem~2.2]{AdCaPi14}, which means that there is a function $f'_{\tilde a}$ such that $c:=f'_{\tilde a}(b)=f_a(b)$ for $b\models p|Ca$. On the other hand, the germ $\tilde a$ lives in $\dcl(Ca)\cap \Int_C(D)$ by Lemma~\ref{L:intgerm} hence we can write $f'_{\theta(a)}$. Then by stable embeddedness, $\tp(c,\theta(a)/C \theta(b))\vdash \tp(c,\theta(a)/C b)$, so $c\in \dcl(\theta(a),\theta(b))$. We can now write $\theta(ab) = c = F(\theta(a),\theta(b))$ for $a\models p$ and $b\models p|Ca$.

By the group chunk theorem (\cite[Proposition~3.15]{HrRi19}), there is a pro-$C$-definable group $\mathfrak g_D$ for which $F$ coincides generically with the multiplication map and $\theta$ extends to a pro-$C$-definable homomorphism $g:G\to \mathfrak g_D$.

Universality follows from the fact that we have taken all of $\dcl(Ca)\cap \Int_C(D)$ in the image of $\theta$.
\end{proof}

\begin{remark}
In the situation of Proposition~\ref{P:metastable-group}, if $C\subseteq C'$, denoting the groups computed using $C$ and $C'$ by $\mathfrak g^{C}_D$ and $\mathfrak g^{C'}_D$, respectively, there is by construction a canonical surjective pro-$C'$-definable homomorphism $\mathfrak g^{C'}_D\rightarrow\mathfrak g^{C}_D$.
\end{remark}

\subsection{Stably dominated groups in $\acvf$}
The following result is a consequence of \cite[Proposition~4.6 \& Corollary~4.11]{HrRi19}.

\begin{fact}\label{F:Finite-index}
Let $N$ be a definable stably dominated group in $\acvf$. Then $N^0=N^{00}$ has finite index in $N$.
\end{fact}

\begin{lemma}\label{L:gen-stdom}
Let $N$ be a definable stably dominated group in $\acvf$. Then the generic types of $N$ are strongly stably dominated.
\end{lemma}
\begin{proof}
By Fact~\ref{F:Finite-index}, we may assume that $N$ is connected.
Let $C$ be a model over which $N$ is defined, and we work over $C$. By \cite[Proposition~4.6]{HrRi19}, there is a definable homomorphism $h:N \to \mathfrak g$, where $\mathfrak g$ is a definable group in $\St_C$ such that the generics of $N$ are stably dominated via $h$, i.e., if $p$ is a generic type of $N$ and $a\models p |C$, then for any 
tuple $b$ with $h(a) \ind_{\St_C(C)} \St_C(b)$, one has $\tp(b/Ch(a))\vdash\tp(b/Ca)$.

We show that $\tp(a/Ch(a))$ is isolated by the formula $h(x)=h(a)$, from which the result follows. By \cite[Lemma~4.9]{HrRi19}, we have that $\tp(a/C)$ is generic in $G$ if and only if $\tp(h(a)/C)$ is generic in $\mathfrak g$. Hence if $h(a')=h(a)$, then $a'$ is generic in $N$ over $C$, hence $\tp(a'/C)=\tp(a/C)$, hence $\tp(a'h(a')/C)=\tp(ah(a)/C)$ and finally $\tp(a'/Ch(a))=\tp(a/Ch(a))$ as required.
\end{proof}

Here is a partial converse to the preceding lemma.

\begin{lemma}\label{L:stdom-def}
Let $N$ be an $\infty$-definable strongly stably dominated connected subgroup of an algebraic group $G$ in $\acvf$. Then $N$ is definable.
\end{lemma}

\begin{proof}
We may replace $G$ by the Zariski closure of $N$ in $G$ and thus assume without loss of generality that $N$ is Zariski dense in $G$. As the generic type of $N$ is strongly stably 
dominated, definability of $N$ follows from \cite[Proposition~4.6 in combination with Corollary~4.16]{HrRi19}.
\end{proof}

\begin{lemma}\label{L:divisible}
Let $N$ be a definable stably dominated group in $\acvf_{0,0}$. Assume that $N$ is abelian and connected. Then $N$ is divisible.
\end{lemma}

\begin{proof}
Let $p$ be a prime number and assume that $pN \neq N$. Then $pN$ has infinite index in $N$ and the group $N/pN$ is an infinite stably dominated abelian $p$-group. We know, e.g., by \cite[Proposition~4.6]{HrRi19}, that this group maps to an infinite stable group ${\mathfrak g}$, where ${\mathfrak g}$ is then also an abelian $p$-group. As ${\mathfrak g}$ is internal to $k$, it is definably isomorphic to a group definable in $ACF_0$, whence definably isomorphic to an algebraic group. It is well known that 
there are no infinite commutative algebraic groups in characteristic 0 which are $p$-groups.\footnote{One may first reduce the statement to linear algebraic groups. For these, by the Jordan decomposition, the result follows from the unipotent and the semisimple cases, which are clear.}
\end{proof}





If $V$ is an algebraic variety defined over a valued field $K$, there is a natural notion of a \emph{bounded} subset of $V(K)$. In \cite[pp.\ 81--83]{Ser90}, this is presented for complete valued fields with archimedean value group, but it readily generalizes to 
arbitrary valued fields, and the results we will use from \cite{Ser90} hold in the general case as well. If $V$ is an affine variety, a subset $B\subseteq V(K)$ is said to be bounded in $V(K)$ if for any $f\in K[V]$ there is $\gamma\in\Gamma_K$ such that 
$\val(f(B))\subseteq[\gamma,\infty]$. If $V$ is an arbitrary variety and $U_1,\ldots,U_m$ an open affine covering of $V$, a subset $B\subseteq V(K)$ is called bounded in $V(K)$ if for $i=1,\ldots,m$ there is a set $B_i\subseteq U_i(K)$ which is bounded in $U_i(K)$ such that $B=\bigcup_{i}B_i$.  We will use the following properties of bounded sets:

\begin{fact}[{\cite[pp.\ 81--82]{Ser90}}]\label{F:Serre-bounded}\mbox{}
\begin{enumerate}
\item If $f:V\rightarrow W$ is a morphism of algebraic varieties and $B\subseteq V(K)$ is bounded in $V(K)$, then $f(B)$ is bounded in $W(K)$.
\item If $V$ is a closed subvariety of $W$ and $B\subseteq V(K)$, then $B$ is bounded in $V(K)$ if and only if $B$ is bounded in $W(K)$.
\item If $f:V\rightarrow W$ is a proper morphism of quasi-projective algebraic varieties and $B\subseteq W(K)$ is bounded in $W(K)$, then $f^{-1}(B)$ is bounded in $V(K)$.\qed
\end{enumerate}
\end{fact}

In what follows, as in \cite{HrLo15}, boundedness  of a definable subset $D$ of $V$ in an algebraic variety $V$ just 
means boundedness of $D(K)$ in $V(K)$, where $K\models\acvf$ and $D$ and $V$ are defined over $K$. This does not depend on the model $K$.

\begin{lemma}\label{L:Bounded-gp-ext}
Let $G$ be an algebraic group, $G_1$ a normal algebraic subgroup and $G_2:=G/G_1$, with projection map $\pi:G\rightarrow G_2$. For a definable subgroup $H$ of $G$, the following are equivalent:
\begin{enumerate}
\item $H$ is bounded in $G$;
\item $H\cap G_1$ is bounded in $G_1$ and $\pi(H)$ is bounded in $G_2$.
\end{enumerate}
\end{lemma}

\begin{proof}
(1)$\Rightarrow$(2) follows from Fact \ref{F:Serre-bounded}(1+2). To prove the converse, we may assume that $H$ is Zariski dense in $G$ and that $G$ is a connected algebraic group. 

\begin{claim*}
There is a definable subset $W$ of $H$ which is bounded in $G$ and such that $\pi(W)=\pi(H)$.
\end{claim*}

To settle the claim, we first note that by the proof of \cite[Corollary~4.5]{HrRi19} there is a definable subset $Y\subseteq H$ such that $\pi(Y)=\pi(H)$ and $\pi\!\!\upharpoonright_Y$ has finite fibers. As algebraic closure (in the field sort) in $\acvf$ is given by field theoretic algebraic closure, we may assume that $Y=H\cap\widetilde{Y}$ for some constructible subset $\widetilde{Y}$ of $G$ such that $\pi\!\!\upharpoonright_{\widetilde{Y}}$ has finite fibers. Since $\pi(\widetilde{Y})$ is a Zariski dense constructible subset of $G_2$, there is an open affine subvariety $U\subseteq G_2$ contained in $\pi(\widetilde{Y})$ such that, setting $Z_U:=\pi^{-1}(U)\cap\widetilde{Y}$, the map $\pi:Z_U\rightarrow U$ is a (surjective and) finite morphism, so is in particular proper. Moreover, note that $\pi(Z_U\cap H)=U\cap\pi(H)$.

As $H$ is Zariski dense in $G$, we may choose finitely many elements $h_1,\ldots,h_n\in H$ such that $(U_i)_{1\leq i\leq n}$ is an (affine) open cover of $G_2$, where $U_i:=\pi(h_i)U$. Thus $\pi(h_iZ_U\cap H)=U_i\cap\pi(H)$ for all $i$. Since $\pi(H)$ is bounded in $G_2$ by assumption, there are bounded definable subsets $B_i$ of $U_i$ such that $\pi(H)=\bigcup_{i=1}^n B_i$. As $\pi:h_iZ_U\rightarrow U_i$ is proper, $W_i:=\pi^{-1}(B_i)\cap h_i Z_U$ is bounded in $h_i Z_U$ by 
Fact \ref{F:Serre-bounded}(3), so bounded in $G$ as well, by Fact \ref{F:Serre-bounded}(2). Thus $W:=\bigcup_{i=1}^n W_i$ is as required, proving the claim.

The set $(H\cap G_1)\times W$ is bounded in $G\times G$, and so $(H\cap G_1)\cdot W$ is  bounded in $G$, as it is the image of a bounded set under  the multiplication map $m:G^2\rightarrow G$. Since $H\subseteq (H\cap G_1)\cdot W$, boundedness of $H$ in $G$ follows.
\end{proof}

\begin{proposition}\label{P:Gdefcompact}
Let $N$ be a definable stably dominated subgroup of an algebraic group $G$. Then $\widehat{N}$ is definably compact.
\end{proposition}

\begin{proof}
We may assume that $G$ is a connected algebraic group. We will first show that $N$ is a bounded subset of $G$. By Chevalley's structure theorem, there is a surjective homomorphism of algebraic groups $\pi:G\twoheadrightarrow A$, with $A$ an abelian variety and $L:=\ker(\pi)$ a linear algebraic group. By \cite[Corollary~4.5]{HrRi19}, $N\cap L$ is a stably dominated group, and so it follows from \cite[Proposition~6.9]{HrRi19} that $(N\cap L)^0$ is a bounded subset of its Zariski closure, which is an algebraic subgroup of $L$. Thus, $(N\cap L)^0$ and so also $N\cap L$ is a bounded subset of $L$. Moreover, $\pi(N)$ is bounded in $A$, as $A$ is a projective variety and so bounded in itself. Thus  $N$ is bounded in $G$ by Lemma~\ref{L:Bounded-gp-ext}, and $N^2$ is bounded in $G^2$.

By Fact~\ref{F:Finite-index}, to prove that $\widehat{N}$ is definably compact, we may assume that $N=N^0$. Let $p_N$ be the unique generic type of $N$. Then the pro-definable set 
$$S=\{r(x,y)\in \widehat{G^2} \mid \widehat{\pi_i}(r)=p_N, i=1,2\}$$ is closed in $\widehat{G^2}$. Moreover, $S\subseteq\widehat{N^2}$, so $S$ is a bounded subset of $\widehat{G^2}$, in the terminology of \cite{HrLo15}. Thus, $S$ is definably compact 
by \cite[Theorem~4.2.19]{HrLo15}. Let $m$ denote the multiplication in $G$. Then $Z=\widehat{m}(S)\subseteq\widehat{G}$ is a definably compact 
subset of $\widehat{N}$ by \cite[Proposition~4.2.9]{HrLo15}. The set $Z$ contains all realized types in $N$ ({\it i.e.}, simple points of $\widehat{N}$), since every element in $N$ is a product of two generics. As the simple points are dense in $\widehat{N}$, we conclude that $\widehat{N}=Z$ is definably compact.
\end{proof}

\subsection{Continuity of the tensor product}\label{Sub:ContTensor}
In this subsection, we continue to work in $\acvf$.

\begin{lemma}\label{L:tensor-cont}
For any definable sets $V$, $W$, the map $\otimes:\widehat{V}\times\widehat{W}\rightarrow\widehat{V\times W}$ is pro-definable.
\end{lemma}
\begin{proof}
Let $p(x)\in \widehat{V}$, $q(y)\in \widehat{W}$ and $r(x,y)=p(x) \otimes q(y)$. Fix some formula $\phi(x,y;t)$. Write $d_p \phi(x; y, t) = \theta_b(y,t)$, where $\theta$ depends only on $\phi$ and $b$ is the canonical parameter of $\theta_b$, so that $b$ is definable from $p$ seen as a point in $\widehat{V}$. For any parameter $c$, we have $\phi(x,y;c)\in r \iff c\models d_q \theta_b(y;t)$. This is definable from $(p,q)\in \widehat{V}\times\widehat{W}$.
\end{proof}

Assume $H$ is a finite-dimensional vector space over $K\models\acvf$. A \emph{semi-lattice} in $H$ is an $\O$-submodule $u\leq H$ for which there is a vector subspace $U\leq H$ such that 
$u/U$ is a lattice in $H/U$. Clearly, the set of semi-lattices $L(H)$ is definable. The \emph{linear topology} on $L(H)$ is the definable topology whose pre-basic open sets are of the form $\Omega_h=\{u\in L(H)\mid h\not\in u\}$ or 
of the form $\Theta_h=\{u\in L(H)\mid h\in\m u\}$. 

We leave the easy proof of the following lemma to the reader. (See \cite[Chapter~5]{HrLo15} for details concerning semi-lattices.).
\begin{lemma}\label{L:properties-semilattices}
\begin{enumerate}
\item Any semi-lattice $u\in L(H)$ is isomorphic to $\O^l\times K^r$, where $l+r=n:=\dim(H)$. In particular, there is a basis $\overline{b}$ of $H$ such that $u$ is diagonal for $\overline{b}$, i.e., $u=\bigoplus_{i=1}^n(u\cap Kb_i)$.
\item If $H=\bigoplus_{i=1}^sH_i$, then $\prod_{i=1}^s L(H_i)$ embeds naturally into $L(H)$. This embedding is a definable homeomorphism onto a closed subset of $L(H)$.
\item Let $\overline{b}=(b_1,\ldots,b_n)$ be a basis of $H$. Then the map $d:\Gamma_\infty^n\rightarrow L(H)$, 
\[d(\gamma_1,\ldots,\gamma_n):=\left\{\sum_{i=1}^n c_ib_i\mid \val(c_i)+\gamma_i\geq0 \text{ for $i=1,\ldots,n$}\right\},\] is a definable homeomorphism 
onto the closed subset $\Delta_{\overline{b}}$ of $L(H)$ consisting of all semi-lattices in $H$ which are diagonal for  $\overline{b}$. In particular, $\Delta_{\overline{b}}$ is a $\Gamma$-internal definable subset of $L(H)$.\qedhere
\end{enumerate}
\end{lemma}

The following strong result about simultaneous diagonalisation is a converse to part (3) of the preceding lemma.

\begin{fact}[\mbox{\cite[Lemma 6.2.2]{HrLo15}}]\label{F:Diagonalization}
Let $X\subseteq L(H)$ be a definable $\Gamma$-internal subset. Then there are finitely many $K$-bases $\overline{b}^1,\ldots,\overline{b}^N$ of $H$ such that any $u\in X$ is diagonal for some $\overline{b}^{i}$.
\end{fact}

\begin{lemma}\label{L:tensor-semi-lattice}
Let $H$ and $H'$ be finite-dimensional vector spaces.
\begin{enumerate}
\item If $u\in L(H)$ and $u'\in L(H')$, then $u\otimes u'\in L(H\otimes H')$, where the tensor product is taken over $\O$. Furthermore, the map 
$\tau:L(H)\times L(H')\rightarrow L(H\otimes H')$, $(u,u')\mapsto u\otimes u'$, is definable.
\item Let $\Sigma\subseteq L(H)$ be a $\Gamma$-internal definable subset. Then the restriction of $\tau$ to $\Sigma\times L(H')$ is continuous 
in the linear topology.
\end{enumerate}
\end{lemma}

\begin{proof}
(1) Clear.

(2) Let $\overline{b}=(b_1,\ldots,b_n)$ be a basis of $H$. We first show the continuity result in the special case where $\Sigma=\Delta_{\overline{b}}\subseteq L(H)$. Let $H_i=Kb_i$ and $E_i=H_i\otimes H'\leq H\otimes H'$. It follows from Lemma~\ref{L:properties-semilattices} that 
$\Delta_{\overline{b}}=\prod_{i=1}^n L(H_i)\cong\Gamma_{\infty}^n$ canonically, and that the map $\tau\upharpoonright_{\Delta_{\overline{b}}\times L(H')}$ decomposes topologically into the product of maps $\tau_i:L(H_i)\times L(H')\rightarrow L(E_i)$. The continuity of $\tau_i$ may 
be easily checked directly. We leave the verification to the reader. This shows the result for $\Sigma=\Delta_{\overline{b}}$.

To prove the general case, using Fact~\ref{F:Diagonalization}, we may find bases $\overline{b}^1,\ldots,\overline{b}^N$ of $H$ such that $\Sigma\subseteq\bigcup_{i=1}^N\Delta_{\overline{b}^{i}}$. Let $(u,u')\in\Sigma\times L(H')$, and let $U$ be an open neighborhood of 
$u\otimes u'$ in $L(H\otimes H')$. Let $I$ be the set of indices $i\in\{1,\ldots,N\}$ such that $u\in\Delta_{\overline{b}^{i}}$. For $i\in I$, choose open neighborhoods $\Omega_i$ of $u$ and $\Omega_i'$ of $u'$ such that $\tau((\Omega_i\cap\Delta_{\overline{b}^{i}})\times\Omega_i')\subseteq U$. For $i\notin I$, we set $\Omega_i=L(H)\setminus\Delta_{\overline{b}^{i}}$ 
(which is an open neighborhood of $u$ by Lemma~\ref{L:properties-semilattices}) and $\Omega_i'=L(H')$. Let $\Omega=\bigcap_{i=1}^N\Omega_i$ and $\Omega'=\bigcap_{i=1}^N\Omega_i'$. By construction, we get $\tau((\Omega\cap\Sigma)\times\Omega')\subseteq U$.
\end{proof}

For $n,d\geq0$ denote by $H_{n,d}$ the vector space of all polynomials of degree $\leq d$ in $n$ variables over the valued field. In what follows, $L(H_{n,d})$ will be endowed with the linear topology. For $p\in \widehat{\mathbb{A}^n}$ let $J_{n,d}(p):=\{h\in H_{n,d}\mid \val(h(p))\geq0\}$. Then 
$J_{n,d}(p)$ is a definable $\O$-submodule of $H_{n,d}$, and it is easy to see that it is a semi-lattice. 

\begin{fact}[\mbox{\cite[Theorem 5.1.4]{HrLo15}}]\label{F:semilattice}
The maps $J_{n,d}:\widehat{\mathbb{A}^n}\rightarrow L(H_{n,d})$ are pro-definable and continuous, and the induced map $J_n:=(J_{n,d})_{d}:\widehat{\mathbb{A}^n}\rightarrow\varprojlim_{d}L(H_{n,d})$ is a pro-definable homeomorphism onto its image.
\end{fact}

\begin{proposition}\label{P:continuity-tensor}
Let $V, W$ be 
algebraic varieties and $\Sigma\subseteq\widehat{V}$ an iso-definable $\Gamma$-internal subset. 
Then the map $\otimes:\Sigma\times\widehat{W}\rightarrow\widehat{V\times W}$ is a (pro-definable) homeomorphism onto its image.
\end{proposition}

\begin{proof}Clearly, the map $\otimes:\Sigma\times\widehat{W}\rightarrow\widehat{V\times W}$ is injective and has a continuous inverse, since the inverse is given by the restriction of the canonical map $\widehat{V\times W}\rightarrow\widehat{V}\times\widehat{W}$. 
It thus suffices to show that $\otimes:\Sigma\times\widehat{W}\rightarrow\widehat{V\times W}$ is continuous. Passing to affine charts and some ambient affine spaces, we may assume that $V=\mathbb{A}^n$ and $W=\mathbb{A}^m$.  

By Fact~\ref{F:semilattice}, it suffices to show that the map $f_d=J_{n+m,d}\circ\otimes:\widehat{\mathbb{A}^n}\times\widehat{\mathbb{A}^m}\rightarrow L(H_{n+m,d})$ restricted to $\Sigma\times\widehat{\mathbb{A}^m}$ is continuous for every $d$. It follows from the definitions of the tensor product of generically stable types and of the maps $J_{l,d}$ 
that $f_d$ factors through $L(H_{n,d})\times L(H_{m,d})$. More precisely, $f_d$ decomposes as 
\[\widehat{\mathbb{A}^n}\times\widehat{\mathbb{A}^m}\xrightarrow{J_{n,d}\times J_{m,d}}L(H_{n,d})\times L(H_{m,d})\xrightarrow{\tau}L(H_{n,d}\otimes H_{m,d})\xrightarrow{\rho}L(H_{n+m,d}),\] 
where $\tau$ is as in Lemma~\ref{L:tensor-semi-lattice} and $\rho$ is the (continuous definable) map induced by the inclusion $H_{n+m,d}\subseteq H_{n,d}\otimes H_{m,d}$.

As $J_{n,d}(\Sigma)$ is a $\Gamma$-internal definable subset of $L(H_{n,d})$, the continuity of the map $f_d\upharpoonright_{\Sigma\times\widehat{\mathbb{A}^m}}$ follows from Lemma~\ref{L:tensor-semi-lattice}.
\end{proof}

\section{Existence of a definable equivariant retraction}\label{S:AllChar}
For our construction of an equivariant deformation retraction, we will only need a special case of the main result of \cite{HrLo15}, namely Fact~\ref{F:Smooth-retract} below, which is much easier to prove than the general non-smooth case. Moreover, all deformation retractions which will appear in our paper are of a very special form. In order to 
simplify the exposition, we give them a name.

\begin{definition}\label{D:SpecialRetr}
Let $V$ be an 
algebraic variety defined over the valued field $F$, and let $X\subseteq V$ be an $F$-definable subset. A pro-$F$-definable continuous map 
$$H:[0,\infty]\times\widehat{X}\rightarrow \widehat{X}$$
is called an \emph{$F$-definable special deformation retraction of $\widehat{X}$ with final image $\Sigma_0$} if the following properties hold:
\begin{enumerate}[(i)]
\item $H_{\infty}=\mathrm{id}_{\widehat{X}}$\label{H-i}
\item $H_0(\widehat{X})\subseteq\Sigma_0$\label{H-ii}
\item $H_t\!\!\upharpoonright_{\Sigma_0}=\mathrm{id}_{\Sigma_0}$ for all $t\in[0,\infty]$.\label{H-iii}
\item For every open subvariety $U$ of $V$, the set $\widehat{U}\cap\widehat{X}$ is invariant under $H$.\label{H-iv}
\item $H(0,x)=H(0,H(t,x))$ for any $x\in\widehat{X}$ and any $t\in[0,\infty]$.\label{H-v}
\item $X^{\#}$ is invariant under $H$.\label{H-vi}
\item $H_0(X)=H_0(\widehat{X})=\Sigma_0$\label{H-vii}
\item $\Sigma_0$ is $F$-definably homeomorphic to a definable subset of $\Gamma^w$, for some finite $F$-definable set $w$.\label{H-viii}
\item For any $x\in \widehat{X}$ and any $t<\infty$, $H(t,x)$ is Zariski generic in $X$.\label{H-ix}
\setcounter{saveenum}{\value{enumi}}
\end{enumerate}
\end{definition}

Here, we denoted by $H_t$ the map sending $x$ to $H(t,x)$.

\begin{remark*}
Let $H$ be  an $F$-definable special deformation retraction of $\widehat{X}$ with final image $\Sigma_0$. Then (\ref{H-vi}) combined with (\ref{H-vii}) and (\ref{H-ix}) yield the following:
\begin{enumerate}[(i)]
\setcounter{enumi}{\value{saveenum}}
\item $\Sigma_0\subseteq X^{\#}$, and every element of $\Sigma_0$ is Zariski generic in $X$.\label{H-x}
\end{enumerate}
\end{remark*}

\begin{fact}\label{F:Smooth-retract}
Let $V$ be a smooth irreducible quasi-projective variety over a valued field $F$ and let $X\subseteq V$ be a $v$-clopen $F$-definable subset such that $\widehat{X}$ is definably compact and definably connected.

Then there is an $F$-definable special deformation retraction $H:[0,\infty]\times \widehat{X}\rightarrow\widehat{X}$ with final image $\Sigma_0$.
\end{fact}

\begin{proof}
That such a map exists with properties (\ref{H-i}-\ref{H-vii}) follows from \cite[Theorem 12.1.1 \& Remark 12.2.3]{HrLo15}. Moreover, as $\widehat{X}$ is definably connected, one may achieve in addition property (\ref{H-viii}), by the discussion in \cite[top of p.~178]{HrLo15}. 

Finally, let us argue why one may achieve property (\ref{H-ix}). It follows from the proof of \cite[Theorem 12.1.1]{HrLo15} that for any $x\in \widehat{X}$ and any $t<\infty$,  the type $H(t,x)$ is Zariski generic in $V$, as this is true for the inflation homotopy $H_{inf}$ defined in \cite[p.~180]{HrLo15} and $H$ is constructed as the composition of $H_{inf}$ with finitely many additional homotopies which are all Zariski generalizing. 
\end{proof}

\begin{lemma}\label{L:Zargen}
Let $V$ be an irreducible variety defined over a model of $\acvf$, and let $\eta:[0,\infty]\rightarrow\widehat{V}$ be a continuous definable path with $\eta(\infty)$ a Zariski generic type in $V$. 
Then $\eta$ is constant in a neighborhood of $\infty$.
\end{lemma}

\begin{proof}
Let $p=\eta(\infty)$. For any open subvariety $U$ of $V$, there exists $\gamma_U<\infty$ such that $\eta([\gamma_U,\infty])\subseteq\widehat{U}$. By the construction of $\widehat{V}$, it is enough to show that for any affine open subvariety $U$ and any regular function $f$ on $U$, the continuous definable function $T:[\gamma_U,\infty]\rightarrow\Gamma_{\infty}$, $T(t):=\val(f_*(\eta(t)))$ 
is constant in a neighborhood of $\infty$. If $f$ is not identically 0 on $U$, then $\val(f_*(p))\neq\infty$, since $p$ is Zariski generic in $U$ by assumption. The result now follows from quantifier elimination in $\Gamma_{\infty}$.
\end{proof}

\begin{corollary}\label{C:paths-smooth}
Let $V$ be a smooth irreducible quasi-projective variety defined over a valued field $F$ and let $X\subseteq V$ be a $v$-clopen $F$-definable subset such that $\widehat{X}$ is definably compact and definably connected. Let $q_0,q_\infty\in\widehat{X}(F)$ with $q_0$ Zariski generic in $V$. Then there is an $F$-definable continuous path $\eta:[0,\infty]\rightarrow\widehat{X}$ with $\eta(0)=q_0$, $\eta(\infty)=q_\infty$ and $\eta(t)$ Zariski generic in $V$ for any $t<\infty$. 

Moreover, if $q_0,q_\infty\in X^{\#}$, then $\eta$ may be required to have its image in $X^{\#}$.
\end{corollary}

\begin{proof}
Choose $H:[0,\infty]\times \widehat{X}\rightarrow\widehat{X}$ as in Fact~\ref{F:Smooth-retract}. Let $r_t = H(t,q_0)$ and $u_t=H(t,q_\infty)$. By Lemma~\ref{L:Zargen}, there is $\gamma<\infty$ such that $r_t=r_\gamma$ for all $t\geq\gamma$. Moreover, since $\widehat X$ is definably connected and definably compact, so is $\Sigma_0=H(0,\widehat{X})$. As $\Sigma_0$ is definably homeomorphic to a subspace of $\Gamma^w$ for some finite set $w$, $\Sigma_0$ is definably path-connected and in particular there is a continuous definable path $s:[\alpha,\beta]\rightarrow \Sigma_0$ between $r_0$ and $u_0$, for some $\alpha<\beta<\infty$, i.e., parametrized by an interval of finite length. We may glue these three paths, which yields a definable path 
$\eta:[0,\infty]\rightarrow\widehat{X}$ between $q_0$ and $q_\infty$ as required. (Zariski genericity of $\eta(t)$ for $t<\infty$ is a consequence of the fact that $H$ satisfies property~(\ref{H-ix}) in Definition~\ref{D:SpecialRetr}.)

The moreover part follows, as $H$ preserves $X^{\#}$ and $\Sigma_0\subseteq X^{\#}$.
\end{proof}

\begin{lemma}\label{L:gp-connected}
Let $G$ be an algebraic group defined over a model of $\acvf$, and let $H=H^0$ be a connected definable subgroup of $G$. Then $\widehat H$ is definably connected.
\end{lemma}

\begin{proof}
The group $H$ acts transitively and definably on the (finite) set of definable connected components of $\widehat H$. The result follows, looking at stabilizers.
\end{proof}

\begin{theorem}\label{P:Group-path-general}
Let $G$ be an algebraic group defined over a valued field $F$, and let $N$ be a commutative $F$-definable stably dominated connected subgroup of $G$. Let $p_e$ be the type of the identity in $N$, and let $p_N$ be the generic type of $N$. Then there is an $F$-definable path 
$q:[0,\infty]\rightarrow N^{\#}$ such that 
\begin{enumerate}
\item [(i)] $q_{\infty}=p_e$ and $q_{0}=p_N$; 
\item[(ii)] for any $t<\infty$, $q_t$ is the generic type of a connected stably dominated definable subgroup $N_t$ of $N$ which is Zariski dense in $N$.
\end{enumerate}
\end{theorem}

\begin{proof}
First note that $p_e$ is strongly stably dominated, as is $p_N$, by Lemma~\ref{L:gen-stdom}, so $p_e,p_N\in N^{\#}(F)$. Furthermore, we may assume that $G$ equals the Zariski closure of $N$, so in particular $G$ is connected and $N$ is Zariski dense in $G$, hence $N$ is $v$-clopen in $G$, as it is a subgroup with non-empty interior for the valuation topology. Moreover, $\widehat{N}$ is definably compact by Proposition~\ref{P:Gdefcompact} and definably connected by Lemma~\ref{L:gp-connected}. We may thus infer from Corollary~\ref{C:paths-smooth} that there is an $F$-definable path  $\eta:[0,\infty]\rightarrow N^{\#}$, such that $\eta_{\infty}=p_e$ and $\eta_{0}=p_N$, with $\eta_t$ strongly stably dominated and Zariski generic in $N$ for any $t<\infty$.

For $p,q\in \widehat{G}$ we denote by $p^{-1}\in\widehat{G}$ the image of $p$ under the inversion map in $G$, and we set $p\ast q:=\widehat{m}(p\otimes q)\in\widehat{G}$, where $m$ denotes the multiplication map in $G$. Now fix $n$ large enough and define $q_t = \eta_t\ast \eta_t^{-1} \ast \cdots$ where the product contains $n$ instances of $\eta_t\ast \eta_t^{-1}$. Then  the map $t\mapsto q_t$ is continuous by Proposition~\ref{P:continuity-tensor}. Moreover, it follows from \cite[Proposition~2.6.12]{HrLo15} that $q_t$ is strongly stably dominated. By \cite[Lemma~5.1]{HrRi19}, $q_t$ is the generic of an $\infty$-definable connected subgroup $N_t$ of $N$. By Lemma~\ref{L:stdom-def}, $N_t$ is in fact definable.  Clearly, $q_\infty=p_e$ and $q_0=p_N$, since $p_e$ and $p_N$ are generic types of connected stably dominated definable subgroups. Moreover, the type $q_t$, obtained as an alternating sum of $\eta_t$, is strongly stably dominated and Zariski generic in $N$ for any $t<\infty$.
\end{proof}




\begin{theorem}\label{T:Contraction-gp}
Let $G$ be an algebraic group defined over a valued field $F$, and let $N$ be a commutative $F$-definable stably dominated connected subgroup of $G$.

Then there is an $F$-definable special deformation retraction $\rho: [0,\infty] \times \widehat N \to \widehat N$ with final image $\{p_N\}$ such that $\rho$ is equivariant under the action of $N$ by multiplication and for each $t<\infty$, $q_t=\rho(t,p_e)$ is the generic type of a connected strongly stably dominated definable subgroup $N_t$ of $N$ which is Zariski dense in $N$. 
\end{theorem}

\begin{proof}
Let $q:[0,\infty]\rightarrow N^{\#}$ be a definable path  as in the conclusion of Theorem~\ref{P:Group-path-general}. Define $r:[0,\infty]\times \widehat{N}\rightarrow\widehat{N\times N}$, $(t,a)\mapsto q_t\otimes a$. Then $r$ is continuous by Proposition~\ref{P:continuity-tensor}, and so is $\rho:=\widehat{m}\circ r:[0,\infty]\times \widehat{N}\rightarrow\widehat{N}$. Clearly, $\rho$ is $N$-equivariant with final image $\{p_N\}$ and satisfies all the required properties from Definition~\ref{D:SpecialRetr}.
\end{proof}



For the definition of the \emph{limit stably dominated subgroup} of a pro-definable group we refer to \cite[Definition~5.6]{HrRi19}.

\begin{lemma}\label{L:towers-limit-stdom}
Let $G$ be an algebraic group, $G_1$ a normal algebraic subgroup and $G_2:=G/G_1$, with projection map $\pi:G\rightarrow G_2$. Suppose that in $G$  the limit stably dominated subgroup exists. Let $N\leq G$ be this group, and suppose that $G/N$ is $\Gamma$-internal. Set $N_1:=N\cap G$ and $N_2:=\pi(N)$. Then the following holds:
\begin{enumerate}
\item For $i=1,2$, $N_i$ is the limit stably dominated subgroup of $G_i$, and $G_i/N_i$ is $\Gamma$-internal.
\item If $N_1$ and $N_2$ are stably dominated, so is $N$.
\end{enumerate}
\end{lemma}

\begin{proof}
(1) Let $N=\bigcup_{t\models q} S_t$, where $(S_t)_{t\models q}$ is a limit stably dominated family for $G$ in the sense of \cite[Definition~5.6]{HrRi19}. Then $S_t\cap G_1$ is stably dominated for any $t\models q$, by \cite[Corollary~4.5]{HrRi19}. Moreover, any connected stably dominated subgroup of $G_1$ is a subgroup of $N$. Thus, $N_1$ is the limit stably dominated subgroup of $G_1$. The group $G_1/N_1$ is $\Gamma$-internal, as it embeds into $G/N$.

Clearly, any $\pi(S_t)$ is stably dominated and connected. Moreover, since $G_2/N_2=\pi(G/N)$ is $\Gamma$-internal, any stably dominated connected subgroup $H$ of $G_2$ is necessarily a subgroup of $N_2$, as $HN_2/N_2$ is stably dominated and $\Gamma$-internal. It follows that $N_2=\bigcup_{t\models q}\pi(S_t)$ is the limit stably dominated subgroup of $G_2$.

(2) follows from \cite[Lemma~4.3]{HrRi19}.
\end{proof} 

\begin{fact}\label{F:Decomp}
Let $S$ be a semi-abelian variety defined over a valued field $F\subseteq K\models\acvf$. Then there is an $F$-definable decomposition $$0\rightarrow N\rightarrow S\rightarrow\Lambda\rightarrow0,$$ with $\Lambda$ $\Gamma$-internal and $N$ stably dominated, definable and connected. The group $N$ is the unique  maximal definable stably dominated and connected  subgroup of $S$. Moreover, if $S=A$ is an abelian variety, $\Lambda$ is definably compact.
\end{fact}

\begin{proof}
If  $S=A$ is an abelian variety, the statement is proven in \cite[Corollary~6.19]{HrRi19}. If $S=\mathbb{G}_m^n$, then $\mathbb{G}_m^n(\O)$ is stably dominated and connected with quotient group $\Lambda\cong(\Gamma,+)^n$, yielding the case of a torus.

Now let $S$ be an arbitrary semi-abelian variety. As $S$ is commutative, by \cite[Theorem~5.16]{HrRi19}, $S$  admits a (unique) limit stably dominated definable subgroup $N$, with quotient $S/N$ $\Gamma$-internal. The result then follows from Lemma~\ref{L:towers-limit-stdom}, as $S$ is the extension of an abelian variety by a torus.
\end{proof}
 
\begin{theorem}\label{T:Retract-all-char}
Let $S$ be a semi-abelian variety defined over  $F\subseteq K\models\acvf$, and let $0\rightarrow N\rightarrow S\rightarrow\Lambda\rightarrow0$ be the decomposition from Fact~\ref{F:Decomp}. 

Then there is an $F$-definable special deformation retraction $\rho: [0,\infty] \times \widehat S \to \widehat S$ with final image $\Sigma \subseteq \widehat{S}$ such that $\rho$ is equivariant under the action of $S$ by multiplication and  for each $t<\infty$, $q_t=\rho(t,p_e)$ is the generic type of a connected strongly stably dominated definable subgroup $N_t$ of $N$, with $N_t$ Zariski dense in $S$.

Moreover, the morphism $\pi:S\rightarrow\Lambda$ induces definable bijection between $\Sigma$ and $\Lambda$.
\end{theorem}

\begin{proof}
We proceed as in the proof of Theorem~\ref{T:Contraction-gp}. We first define a map $r:[0,\infty]\times \widehat{S}\rightarrow\widehat{S\times S}$, $(t,a)\mapsto q_t\otimes a$, where $q$ is as in Theorem~\ref{P:Group-path-general}. The map $\rho:=\widehat{m}\circ r:[0,\infty]\times \widehat{S}\rightarrow\widehat{S}$ is continuous and satisfies all the required properties from Definition~\ref{D:SpecialRetr}, and it is clearly $S$-equivariant. The restriction $\rho_S=\rho\upharpoonright_{[0,\infty]\times S}$ has final image equal to $\{a+ p_N\mid a\in S\}$, a set which may be identified with $S/N\cong\Lambda$ and which is thus $\Gamma$-internal. 

It follows from the definitions that $\rho$ is the canonical extension (in the sense of \cite[Section~3.8]{HrLo15}) of $\rho_S$. Thus, the final image of $\rho$ is equal to $\Lambda$ as well, as $\widehat{\Lambda}=\Lambda$.
\end{proof}

\section{An explicit definable equivariant retraction in equicharacteristic 0} \label{S:Retract0}
In this section, we will give an alternative, more explicit, construction of an equivariant  definable special deformation retraction in the case of equicharacteristic 0. This construction does not require knowing in advance that a (non-equivariant) retraction exists. More importantly, it does not require that the stably dominated connected group $N$ be commutative, in order to show that 
its stable completion $\widehat{N}$ allows for an $N$-equivariant definable special deformation retraction to the generic type of $N$. We believe that even in the commutative case, it might be useful in its own right.

\subsection{Internality of quotients in ACVF$_{0,0}$}
The following result will be used to define an intrinsic scale, given by subgroups, in any stably dominated definable subgroup of an algebraic group in a model of $\acvf_{0,0}$.

\begin{lemma}\label{lemma:internal}
Work in $\acvf_{0,0}$. Let $D,D'$ be definable subgroups of $\O$. Then the following are equivalent:
\begin{enumerate}
\item $D'\subseteq D$.
\item $\O/D$ is $\O/D'$-internal.
\item $\O/D$ is almost $\O/D'$-internal.
\end{enumerate}
\end{lemma}

\begin{proof}
$1.\Rightarrow 2.\Rightarrow 3.$ is clear. 

We now prove $3.\Rightarrow1.$ Note that every definable subgroup of $\O$ is of the form $\gamma\m$ or $\gamma\O$ for some $0\leq\gamma\leq\infty$, so in particular 
the set of definable subgroups of $\O$ is totally ordered by inclusion. Whenever $D\subsetneq D'$ are definable subgroups, there is $0<\gamma\in\Gamma$ 
such that $D\subseteq\gamma\m\subsetneq\gamma\O\subseteq D'$. It is  thus enough to show that for any $0<\gamma\in\Gamma$, the 
set $\O/\gamma\m$ is not almost $\O/\gamma\O$-internal. 

The idea of the proof is similar to \cite[Lemma~5.1]{HrTa06}. Consider the field of generalized power series $K=k((X^\Gamma))$, where $k$ is an algebraically closed field of characteristic 0 and $\Gamma$ a non-trivial divisible ordered abelian group. Let $l=k((T^\QQ))$ (with the trivial valuation) and consider $L=l((X^\Gamma))$ with the $X$-adic valuation, an elementary extension of $K\models\acvf_{0,0}$. 

Given $0<\gamma\in\Gamma(K)=\Gamma$ and $\rho\in L$ of valuation 0 (e.g., $\rho\in k^{\times}$) we now define an automorphism 
$\sigma_\rho$ of $L$ fixing $K\cup(\O/\gamma\O)(L)$ pointwise such that for distinct $\rho,\rho'\in k^{\times}$, we have 
$\sigma_\rho(T+\gamma\m)\neq\sigma_{\rho'}(T+\gamma\m)$. This will show the result, as then $(\O/\gamma\m)(L)\subsetneq\acl(K\cup(\O/\gamma\O)(L))$.

Define $\sigma=\sigma_\rho$ on monomials as follows:
$$\sigma(aT^qX^g):=aT^qX^g\exp(q\rho X^\gamma).$$
Then extend the map by linearity to generalized series. We check that $\sigma$ is an automorphism having the required properties:

(1) $\sigma$ is a continuous automorphism: linearity and continuity follow from the definition. Multiplicativity can be checked on monomials: $\sigma(a_1T^{q_1}X^{g_1}) \cdot \sigma(a_2T^{q_2}X^{g_2}) = a_1a_2T^{q_1+q_2}X^{g_1+g_2}\exp(q_1\rho X^{\gamma})\exp(q_2\rho X^{\gamma}) = \sigma(a_1a_2T^{q_1+q_2}X^{g_1+g_2})$.

For any $x\in\O$, we have \begin{equation}\label{eq:sigma}v(x-\sigma(x)) \geq v(x)+\gamma.\end{equation} This, along with continuity, implies that $\sigma$ is a bijection: if we have $y$ such that $v(\sigma(y)-x) =: \eta$, we let $y' = y+x-\sigma(y)$ and obtain a better candidate for a preimage of $x$, that is $v(\sigma(y')-x) > v(\sigma(y)-x)$. By transfinite induction, we find an actual preimage of $x$.

(2) $\sigma$ fixes $K\cup (\O/\gamma\O)(L)$ pointwise: For $x\in K$, $\sigma(x)=x$ is immediate from the definition and \eqref{eq:sigma} gives the rest of the statement.

(3) $\sigma_\rho(T+\gamma\m)\neq\sigma_{\rho'}(T+\gamma\m)$ for $\rho \neq \rho'$ follows  from the computation $v(\sigma_\rho(T)-\sigma_{\rho'}(T)) = v((\rho-\rho') TX^\gamma) = \gamma$.
\end{proof}

\begin{remark}
In positive and in mixed characteristic, the statement of Lemma~\ref{lemma:internal} does not hold: in $\acvf_{p,p}$, given any $\gamma\in\Gamma_{>0}$ and $c\in\O$ with $v(c)=\gamma$, the Frobenius automorphism induces a group isomorphism $\O/\gamma\O=\O/c\O\cong\O/c^p\O=\O/(p\gamma)\O$; in $\acvf_{0,p}$, setting $\gamma_0:=\val(p)\in\Gamma_{>0}$, the map $x\mapsto x^p$ induces a definable surjection  $\O/\gamma_0\O\twoheadrightarrow\O/(2\gamma_0)\O$.
\end{remark}

\begin{corollary}\label{Cor:multi-int}
Work in $\acvf_{0,0}$. Let $D$ be a definable subgroup of $\O$, and let $d\in\NN$. Then $(\O/D)^d$ is the maximal (almost) $\O/D$-internal quotient of $\O^d$.
\end{corollary}

\begin{proof}
It is clear that $(\O/D)^d$ is $\O/D$-internal. Conversely, let $\O^d/N$ be (almost) $\O/D$-internal for some ($\infty$-)definable subgroup $N\leq \O^d$. For 
$i\in\{1,\ldots d\}$, let $B_i=\{0\}^{i-1}\times\O\times\{0\}^{d-i-1}\cong\O$. It follows that for any $i$ the group $B_i/(N\cap B_i)$ is (almost) $\O/D$-internal, and 
so $\{0\}^{i-1}\times D\times\{0\}^{d-i-1}\leq N\cap B_i$ by Lemma~\ref{lemma:internal}. Thus $D^d\leq N$.
\end{proof}

\begin{corollary}\label{Cor:finite-gen}
Work in $\acvf_{0,0}$. Let $C$ be a model  (or more generally a base structure consisting of a field). Let $\gamma\in \Gamma$, and set $C_\gamma:=\acl(C\gamma)$. Let $D$ be a $C_\gamma$-definable subgroup of $\O$. If $a=(a_1,\ldots,a_d)$ is generic in $\O^d$ over $C$, the tuple $(a_1/D,\ldots,a_d/D)$ $\dcl$-generates $\dcl(C_\gamma a)\cap\Int_{C_\gamma}(\O/ D)$ over $C_\gamma$.


More generally, if $a\in K^m$ realizes a strongly stably dominated type over $C$, there is a tuple $b=(b_1,\ldots,b_d)$ from $C(a)$ such that $b$ is generic in $\O^d$ over $C$ and $\acl(C_\gamma a)\cap\Int_{C_\gamma}(\O/D)$ is finitely $\acl$-generated by $(b_1/D,\ldots,b_d/D)$  over $C_\gamma$. 
\end{corollary}

\begin{proof}
First suppose that $a$ is generic in $\O^d$ over $C$. By the previous corollary, $(\O/D)^d$ is the maximal (almost) $\O/D$-internal quotient of $\O^d$, over any set of parameters. As $\tp(a/C)$ is stably dominated and thus orthogonal to $\Gamma$, the tuple $a$ is generic in $\O^d$ over $C_\gamma$. By the last sentence of Proposition~\ref{P:metastable-group}, $\dcl(aC_\gamma)\cap\Int_{C_\gamma}(\O/D)$ 
is then interdefinable with $(a_1/D,\ldots,a_d/D)$ over $C_\gamma$, proving the result.

Now suppose $a\in K^m$ with $\tp(a/C)$ strongly stably dominated, i.e., $\trdeg(a/C)=\trdeg(k_{C(a)}/k_{C})=d$ by Fact~\ref{fact:ssd}. We find $b\in C(a)$ generic 
in $\O^d$ over $C$. As $b$ and $a$ are interalgebraic over $C$, we have $\acl(C\gamma b)\cap\Int_{C\gamma}(\O/D)=\acl(C\gamma a)\cap\Int_{C\gamma}(\O/D)$, so the result follows from the special case.
\end{proof}

\subsection{Linearization}\label{Sub:Linear}
In this subsection, we work in ACVF$_{0,0}$. Let $C$ be a base structure consisting of a field. We fix a $C$-definable stably dominated connected subgroup $N$ of an algebraic group. For $\gamma \in \Gamma$, let $C_\gamma = \acl(C\gamma)$. Set $d:=\dim(N)$.

\begin{fact}[\mbox{\cite[Section 7]{HrTa06}}]\label{F:balls-SE}
For $\gamma\in \Gamma$, $\O/\gamma\O$ and $\O/\gamma \m$ are stably embedded. 
\end{fact}

By Proposition~\ref{P:metastable-group}, we let $N_\gamma$ be the kernel of the maximal $\O/\gamma \O$-internal quotient of $N$ and similarly $N^+_\gamma$ is defined as the kernel of the maximal $\O/\gamma\m$-internal quotient of $N$, both computed over $C_\gamma$.

The groups $N_\gamma$ and $N_\gamma^+$ are $\infty$-definable. Note that if $\gamma<\delta$, then $N_\gamma \supseteq N^+_\gamma \supseteq N_\delta \supseteq N^+_\delta$.

\begin{lemma}\label{L:Ng_def}
The groups $N_\gamma$ and $N^+_\gamma$ have bounded index inside some definable group and are intersections of definable groups.
\end{lemma}

\begin{proof}
By Fact~\ref{F:Generics}(\ref{I:ProGroup}), the quotient $N/N_\gamma$ can be written as an inverse limit of definable groups $\pi_i(N)$. By Lemma~\ref{L:gen-stdom}, the generic type $p_N$ of $N$ is strongly stably dominated and by Corollary~\ref{Cor:finite-gen}, over $C_\gamma$, the set  $\acl(C_\gamma a)\cap \Int_{C_\gamma}(\O/\gamma\O)$ is in the algebraic closure of a finite set in $\O/\gamma\O$. By construction of $N_\gamma$, if $a\models p_N | C$, then $\acl(C_\gamma a)\cap \Int_{C_\gamma}(\O/\gamma\O)$ is interalgebraic over $C_\gamma$ with the sequence $\pi_i(a)$. Hence we can find $i_0$ such that $\acl(C_\gamma a)\cap  \Int_{C_\gamma}(\O/\gamma\O)$ is already interalgebraic with $\pi_{i_0}(a)$ over $C_\gamma$. Then $N_\gamma$ has bounded index inside the kernel of $\pi_i$. Also $N_\gamma$ is the intersection of the kernels of $\pi_i$ which are definable subgroups of $N$.

The same arguments work for $N^+_\gamma$.
\end{proof}

\begin{lemma}\label{L:Ng_metastable}
The group $N_\gamma$ is strongly stably dominated.
\end{lemma}

\begin{proof}
We first show that the quotient $\mathfrak{g}_\gamma:=N_\gamma/N_\gamma^+$ is stable of Morley rank $d$. Let $p=p_N$ be the generic type of $N$, and let $\Lin(p)$ be the generic type of $\O^d$. Then any 
$a\models p | C$ is interalgebraic over $C$ with some $b=(b_1,\ldots,b_d)\models \Lin(p) | C$. Denote by $r_\gamma$ (by $r_{\gamma}^+$, respectively) the generic type 
of the maximal $\O/\gamma \O$-internal (maximal $\O/\gamma\m$-internal, respectively) quotient of $N$ over $C$, with realizations $a_\gamma$ and $a_{\gamma}^+$, images of $a$ via the canonical projection. Let $b_\gamma=(b_1/\gamma\O,\ldots,b_d/\gamma \O)$, and define similarly $b_{\gamma}^+$. We have the following:

\begin{itemize}
\item $\acl(C_\gamma b_\gamma^+)=\acl(C_\gamma a_\gamma^+)$;
\item $\acl(C_\gamma b_\gamma)=\acl(C_\gamma a_\gamma)$ and as $\tp(b_\gamma^+/\acl(C_\gamma b_\gamma))$ is interdefinable with the generic of a $k^d$-torsor, then $\tp(a_\gamma^+/\acl(C_\gamma a_\gamma))$ is interalgebraic with the generic of a $k^d$-torsor. In particular, that type is stable of Morley rank $d$.
\item The generic of $\ker(N/N_\gamma^+\rightarrow N/N_\gamma)\cong\mathfrak{g}_\gamma$ is interdefinable with the non-forking extension of $\tp(a_\gamma^+/C_\gamma a_\gamma)$ by Lemma~\ref{lem:coset-generic}, hence $\mathfrak{g}_\gamma$ is a stable group of Morley rank $d$.
\end{itemize}
Any type $p_\gamma$ of $N_\gamma$ projecting on the generic of $\mathfrak g_\gamma$ is strongly stably dominated. In particular, it is definable over $C_\gamma$. The same is true for any translate of $p_\gamma$, since the generic of $\mathfrak g_\gamma$ is invariant under translation. Hence $p_\gamma$ is a generic type of $N_\gamma$ and $N_\gamma$ is strongly stably dominated.
\end{proof}

\begin{corollary}
The groups $N_\gamma$ and $N^+_\gamma$ are definable.
\end{corollary}

\begin{proof}
By Lemma~\ref{L:Ng_metastable}, $N_\gamma$ is strongly stably dominated, so it is definable by Lemma~\ref{L:stdom-def}.

Now by the proof of Lemma~\ref{L:Ng_metastable}, the quotient $N_\gamma/N^{+}_\gamma$ is $\omega$-stable. By Lemma~\ref{L:Ng_def}, write $N^+_\gamma=\bigcap_i H_i$, where the $H_i$'s are definable subgroups of $N_\gamma$. Since there is no infinite descending chain of definable subgroups of $N_\gamma/N^+_\gamma$, the intersection $\bigcap_i H_i$ is equal to a finite subintersection, hence $N^+_\gamma$ is definable.
\end{proof}

\begin{corollary}\label{C:connected}
The group $N^{00}_\gamma=N^0_\gamma$ is of finite index in $N_\gamma$.
\end{corollary}

\begin{proof}
This is just Fact~\ref{F:Finite-index}, combined with Fact~\ref{F:Generics}(\ref{I:Connected-component}).
\end{proof}

\begin{lemma}\label{L:union}
Let $\gamma \in [0,\infty]$. Set $V_\gamma=\bigcup_{\delta>\gamma}N_\delta$ and $W_\gamma = \bigcap_{\delta<\gamma}N_\delta$. Then the following holds:
\begin{enumerate}
\item $V_\gamma=\bigcup_{\delta>\gamma}N^+_\delta\leq N_\gamma^+$ and $(N^+_\gamma:V_\gamma)<\infty$, so in particular $\left(N^+_\gamma\right)^0=V_\gamma^0$.
\item $N_\gamma \leq W_\gamma=\bigcap_{\delta<\gamma}N^+_\delta$ and $(W_\gamma:N_\gamma)<\infty$, so in particular $N^0_\gamma=W_\gamma^0$.
\end{enumerate}
\end{lemma}

\begin{proof}
We know that $V_\gamma=\bigcup_{\delta>\gamma}N^+_\delta\subseteq N_\gamma^+$ and $N_\gamma \subseteq W_\gamma=\bigcap_{\delta<\gamma}N_\delta = \bigcap_{\delta<\gamma}N^+_\delta$.

Let us show that $N^0_\gamma = W_\gamma^0$. We fix some saturated model $M_0$ and we work over $C_\gamma$. Suppose that the index of $N_\gamma$ in $W_\gamma$ is infinite.  Let $a\in M_0$ realize the generic of $N$ over $C_\gamma$ and let $a_\gamma$ be the image of $a$  in $N/N_\gamma$ under the projection map, as in Lemma~\ref{L:Ng_metastable}, and let $\tilde a_\gamma$ be the image of $a$ in $N/W_\gamma$. We have that $\tilde a_\gamma$ is algebraic over $a_\gamma$ and that $a_\delta$, $\delta<\gamma$ is algebraic over $\tilde a_\gamma$. However, by assumption, $a_\gamma$ is not algebraic over $\tilde a_\gamma$ and there exists an automorphism fixing $C_\gamma \tilde a_\gamma$ along with $a_\delta$, $\delta<\gamma$ and for which $a_\gamma$ has an infinite orbit. We now define $b_\delta$, $\delta\leq \gamma$ as in the proof of Lemma~\ref{L:Ng_metastable}. Then each $b_\delta$ is interalgebraic with $a_\delta$. Therefore $b_\gamma$ has an infinite orbit under $\sigma$ whereas each $b_\delta$, $\delta<\gamma$ has a finite orbit. Since $b_{\delta}\in \dcl(C_\gamma b_{\delta'})$ for $\delta<\delta'$, taking a large enough power of $\sigma$, we may assume that $\sigma$ actually fixes each $b_\delta$, $\delta<\gamma$. This is impossible since any element $d$ of $\mathcal O/\gamma \mathcal O$ is determined by the sequence $(d_\delta : \delta<\gamma)$ where $d_\delta$ is the image of $d$ in $\mathcal O/\delta \mathcal O$.

The fact that $(N^+_\gamma:V_\gamma)<\infty$  is proved in a similar way. Assume that $V_\gamma$ has infinite index in $N_\gamma^+$. Let $\tilde a^+_\gamma$ be the image of $a$ in $N/V_\gamma$, and $a^+_\gamma$ that in $N/N^+_\gamma$. Take as above an automorphism $\sigma$ fixing $C_\gamma a^+_\gamma$ and under which $\tilde a^+_\gamma$ has infinite orbit. This implies that each $a_\delta$, $\delta>\gamma$ has infinite orbit. The same is true for each $b_\delta$, $\delta>\gamma$, however $b^+_\gamma$ has finite orbit under $\sigma$. Taking a power of $\sigma$, we may assume that $b^+_\gamma$ is fixed. As $b^+_\gamma=b+\gamma\m$, using $\sigma(\gamma)=\gamma$, we get $\sigma(b)+\gamma\m=\sigma(b^+_\gamma)$ and so $\delta_0=v(b-\sigma(b))>\gamma$. This implies that $\sigma(b_{\delta_0})=b_{\delta_0}$, which is the desired contradiction.
\end{proof}


\begin{corollary}
For any $\gamma\in[0,\infty]$, one has $\left(N^+_\gamma\right)^0=\left(N^+_\gamma\right)^{00}$.
\end{corollary}

\begin{proof}
This follows from $\bigcup_{\delta>\gamma}N^0_\delta=\left(N^+_\gamma\right)^0$ and Corollary~\ref{C:connected} since an increasing union of strongly connected groups is strongly connected.
\end{proof}

\begin{lemma}\label{L:nogenstable}
Let $H$ be a definable stably dominated group such that $N_\gamma\supseteq H\supseteq \left(N_{\gamma}^+\right)^0$. Then $(N_\gamma: H)<\infty$.
\end{lemma}

\begin{proof}
We work over a model $M$ containing $C_\gamma$. Suppose $\acl(a_\gamma)\subsetneq \acl(a')\subseteq\acl(a^+_\gamma)$, where, $a_\gamma$, $a'$ and $a_\gamma^+$ are the images of $a\models p_N | M$ under the canonical projection of $N$ onto $N/N_\gamma$, $N/H$ and $N/N^+_\gamma$, respectively.

As non-forking extensions/restrictions and translates of (strongly) stably dominated types are (strongly) stably dominated, 
it follows from Lemma~\ref{lem:coset-generic} that $\tp(a/a_\gamma)$ is generic in $a+N_\gamma$ and strongly stably dominated. Similarly, $\tp(a/a')$ is generic in $a+H$ and stably dominated.

As $a_\gamma^++N_\gamma^+\in \St_{a_\gamma}(a)$ is of Morley rank $d$ over $a_\gamma$, it follows that $\St_{a_\gamma}(a)$ is interalgebraic over $a_\gamma$ with 
$a_\gamma^+$. Thus there is some non-algebraic element from $\St_{a_\gamma}(a)$ which is in $\acl(a')$. Note that $\St_{a_\gamma}(a)$ and $\St_{a'}(a)$ are interdefinable over $a'$, since $a'\in\St_{a_\gamma}(a)$ and $a_\gamma\in\dcl(a')\subseteq\dcl(a)$ (cf. \cite[Lemma 3.5 \& Remark 7.9]{HHM2}). 
It follows that $\dim(\St_{a'}(a))<d$, and so $\tp(a/\acl(a'))$ is not strongly stably dominated by Fact~\ref{fact:ssd}. Thus 
the generic type of $H$ is not strongly stably dominated either, contradicting Lemma~\ref{L:gen-stdom}.
\end{proof}

\begin{remark}
Any $\infty$-definable group $H$ such that $N_\gamma\supseteq H\supseteq N_{\gamma}^+$ is definable, since 
$H/N_{\gamma}^+$ is definable, as it is an $\infty$-definable subgroup of an $\omega$-stable group. 
\end{remark}





\subsection{Definability and continuity}\label{Sub:Def-cont}
We keep the notation and conventions from Subsection~\ref{Sub:Linear}.

\begin{lemma}\label{L:Definable-Ngamma}
There are $C$-definable families $\left(\widetilde{N}_\gamma\right)_{\gamma\in[0,\infty]}$ and $\left(\widetilde{N}^+_\gamma\right)_{\gamma\in[0,\infty]}$ of definable sugroups of $N$ such that for any $\gamma$,  $N_\gamma\leq\widetilde{N}_\gamma$ with $(\widetilde{N}_\gamma:N_\gamma)<\infty$ and similarly $N^+_\gamma\leq\widetilde{N}^+_\gamma$ with $(\widetilde{N}^+_\gamma:N^+_\gamma)<\infty$.

In addition, we may choose these families so that $\widetilde{N}^+_\gamma\leq\widetilde{N}_\gamma\leq \widetilde{N}^+_\delta\leq \widetilde{N}_\delta$ whenever $\gamma<\delta$.
\end{lemma}

\begin{proof}
For $\gamma\in[0,\infty]$, there is a formula  $\psi(x,y)$ with parameters from $C$ such that $\psi(x,\gamma)$ defines $N_\gamma$.  Corollary~\ref{Cor:finite-gen} together with Proposition~\ref{P:metastable-group} implies that there are $\overline{a}$ generic in $N/N_\gamma$ over $C_\gamma$ and $\overline{b}$ generic in $(\O/\gamma\O)^d$ over $C_\gamma$ such that $\overline{b}\in\dcl(C_\gamma\overline{a})$ and $\overline{a}\in\acl(C_\gamma\overline{b})$.  By compactness, there is $\theta(y)\in\tp(\gamma/C)$ such that whenever $\models\theta(\delta)$, then $\psi(x,\delta)$ defines a subgroup $N_\delta'$ of $N$ such that there are $\overline{a}$ generic in $N/N'_\delta$ over $C_\gamma$ and $\overline{b}$ generic in $(\O/\gamma\O)^d$ over $C_\gamma$ with $\overline{b}\in\dcl(C_\gamma\overline{a})$ and $\overline{a}\in\acl(C_\gamma\overline{b})$. It follows from Proposition~\ref{P:metastable-group} that $N_\delta\leq N'_\delta$ and $(N'_\delta:N_\delta)<\infty$ for any such $\delta$. By compactness, we obtain a $C$-definable family $\left(\widetilde{N}_\gamma\right)_{\gamma\in[0,\infty]}$ with the required properties.  

 In exactly the same way, one proves the existence of the family $\left(\widetilde{N}^+_\gamma\right)_{\gamma\in[0,\infty]}$.

In order to achieve the additional requirement, it is enough to replace the group $\widetilde{N}_\gamma$ by $\bigcap_{\gamma\leq\delta}\widetilde{N}_\delta\cap\bigcap_{\gamma<\delta} \widetilde{N}^+_\delta$, and similarly $\widetilde{N}^+_\gamma$ by $\bigcap_{\gamma\leq\delta}\widetilde{N}_\delta\cap\bigcap_{\gamma\leq\delta} \widetilde{N}^+_\delta$.
\end{proof}

\begin{lemma}\label{L:nfcp}
The theory $\acvf$ admits elimination of $\exists^\infty$ in imaginary sorts.
\end{lemma}

\begin{proof}
By  \cite[Lemma 2.6.2]{HHM}, if $D$ is a definable (imaginary) set, then either it is $k$-internal, or,  for some $m$, there is a definable surjective map from $D^m$ to an infinite interval of $\Gamma$. The second case is an open condition and implies that $D$ is infinite. In the first case, since we know that $k$ eliminates $\exists^\infty$, we also see that $D$ being infinite can be expressed as a definable condition on the parameters defining $D$.
\end{proof}

\begin{lemma}\label{L:pro-def}
Let $q_\gamma\in\widehat{N}$ be the generic type of 
$N_\gamma^0$. Then the map $q:[0,\infty]\rightarrow\widehat{N}$, $\gamma\mapsto q_\gamma$ is pro-definable. Its image is iso-definable and $\Gamma$-internal.
\end{lemma}

\begin{proof}
Let  $\left(\widetilde{N}_\gamma\right)_{\gamma\in[0,\infty]}$ be a family as given by Lemma~\ref{L:Definable-Ngamma}. Then $N_\gamma^0=\widetilde{N}^0_\gamma$ for all $\gamma$. For every formula $\phi(x,y)$, the set 
$$\mathrm{Gen}(\phi):=\{(\gamma,p)\in[0,\infty]\times\widehat{N}\,\mid\, p \text{ is $\phi$-generic in $\widetilde{N}_\gamma$}\}$$
is relatively definable in $[0,\infty]\times\widehat{N}$. This can be seen as follows. By generic stability, the $\phi$-definitions 
of elements of $\widehat{N}$ are uniform. Now $\widetilde{N}_\gamma$ acts on the $\phi$-definitions of elements of $\widehat{\widetilde{N}_\gamma}$. By Lemma~\ref{L:nfcp}, there is $n_\phi$ such that whenever 
the orbit of a $\phi$-definition under this action is finite, it is of cardinality at most $n_\phi$. It follows that the set $\mathrm{Gen}:=\bigcap_{\phi}\mathrm{Gen}(\phi)\subseteq[0,\infty]\times\widehat{N}$ is pro-definable.

The set $\mathrm{PrGen}:=\{(\gamma,p)\in \mathrm{Gen}\mid p*p=p\}$ is then pro-definable as well, and it is equal to $\mathrm{Graph}(q)$, since for a generic type $p$ of a stably dominated group $H$, one has $p*p=p$ precisely if $p$ is the principal generic of $H$. This proves pro-definability of $q$.

As $q$ is injective, iso-definability of its image follows from compactness and then this image is $\Gamma$-internal by definition.
\end{proof}

\begin{theorem}\label{P:continuity}
The map $q$ is continuous. Thus, $q:[0,\infty]\rightarrow N^{\#}$ is a definable path between $q_{\infty}=p_e$ and $q_0=p_N$ along generic types of strictly increasing strongly stably dominated connected definable subgroups of $N$, with $q_\gamma$ Zariski generic in $N$ for every $\gamma<\infty$.
\end{theorem}

\begin{proof}
By definable compactness of $\widehat{N_\gamma}$ (Proposition~\ref{P:Gdefcompact}), $l^+_\gamma:=\lim_{\delta\models\gamma^+}q_\delta$ exists and is in $\widehat{N_\gamma}$. We claim that $l^+_\gamma=q_\gamma$. To see this, first note that, for any $\delta>\gamma$ and $a\in N_\delta$, we have $a\cdot q_\delta = q_\delta$, hence also $a\cdot l^+_\gamma=l^+_\gamma$, as the map $q\mapsto a\cdot q$ is continuous (e.g., by Proposition~\ref{P:continuity-tensor}).
Moreover, as $q_\delta * q_\delta = q_\delta$ for all $\delta$, we have $l^+_\gamma = l^+_\gamma * l^+_\gamma$ by continuity of $\otimes$ (Proposition~\ref{P:continuity-tensor}). Hence $l^+_\gamma$ is the generic type of a subgroup of $N_\gamma$ containing $\bigcup_{\delta>\gamma}N_\delta=N_\gamma^{+}$ (Lemma~\ref{L:union}). By Lemma~\ref{L:nogenstable}, $l^+_\gamma$ must be the generic of $N_\gamma$, namely $l^+_\gamma=q_\gamma$ as required.

We next show continuity at $\gamma^-$ (including $\gamma=\infty$): As before, by continuity of $\otimes$, we know that $l^-_\gamma:=\lim_{\delta\models\gamma^-}q_\delta$ 
is an idempotent. Moreover, since $q_\gamma * q_\delta = q_\delta$, for every $\delta<\gamma$, this 
holds in the limit as well, {\it i.e.}, $q_\gamma * l^-_\gamma = l^-_\gamma$. As $l^-_\gamma\in\bigcap_{\delta<\gamma}\widehat{N_\delta}=\widehat{W_\gamma}$, it follows from Lemma~\ref{L:union} that $l^-_\gamma=q_\gamma$, since the only idempotent generic type of $W_\gamma$ is the principal generic.
\end{proof}

The proofs  of Theorem~\ref{T:Contraction-gp} and Theorem~\ref{T:Retract-all-char}, respectively, show that Theorem~\ref{P:continuity} yields the following two corollaries. 

\begin{corollary}\label{C:Contract-char-00}
Let $G$ be an algebraic group defined over a model of $\acvf_{0,0}$, and let $N$ be a $C$-definable stably dominated connected subgroup of $G$.

Then there is a $C$-definable special deformation retraction $\rho: [0,\infty] \times \widehat N \to \widehat N$ with final image $\rho(0,\widehat{N})=\{p_N\}$ such that $\rho$ is equivariant under the action of $N$ by multiplication and for each $t<\infty$, $q_t=\rho(t,p_e)$ is the generic type of a definable strongly stably dominated connected Zariski dense subgroup $N_t$ of $N$, with $N_s\supsetneq N_t$ if $s<t$.  
\end{corollary}

\begin{corollary}\label{C:Retract-char-00}
Let $S$ be a semi-abelian variety defined over $F\subseteq K\models\acvf_{0,0}$, and let $0\rightarrow N\rightarrow S\rightarrow\Lambda\rightarrow0$ be the decomposition from Fact~\ref{F:Decomp}. 

Then there is an $F$-definable special deformation retraction $\rho: [0,\infty] \times \widehat S \to \widehat S$ with final image $\Sigma:=\rho(0,\widehat{S})\subseteq \widehat{S}$ such that $\rho$ is equivariant under the action of $S$ by multiplication and for each $t<\infty$, $q_t=\rho(t,p_e)$ is the generic type of a definable strongly stably dominated connected subgroup $N_t$ of $N$ which is Zariski dense in $S$ and such that $N_s\supsetneq N_t$ whenever $s<t$.  Moreover, $\Sigma$ is in definable bijection with $\Lambda$, canonically.
\end{corollary}

\section{Application to the topology of $S^{\an}$}\label{S:Berkovich}
Let $S$ be a semi-abelian variety defined over a valued field $F$ with $\Gamma_F\leq\RR$. In this section, using  the methods and results from \cite[Chapter 14]{HrLo15}, we will deduce from our results on $\widehat{S}$ the existence of an equivariant strong deformation retraction of (the underlying topological space of) the Berkovich analytification $S^{\an}$ onto its skeleton, which is a simplicial complex carrying the structure of a piecewise linear group.

Recall that if $B$ is a parameter set in some model of $\acvf$, a type $p\in S(B)$ is  \emph{almost orthogonal to $\Gamma$}, denoted by $p\perp^{a}\Gamma$, if $\Gamma(Ba)=\Gamma(B)$ for  $a\models p$. If this is the case then for a $B$-definable function $g:X\rightarrow\Gamma_\infty$, with $p$ concentrating on the definable set $X$, we set $g(p):=g(a)$, where $a\models p$. This is well-defined, as the value $g(a)$ does not depend on the realization $a$.

Until the end of this section, we fix a valued field $F$  with $\Gamma_F\leq\RR$. We will now introduce some objects and notation (mostly) from \cite[Chapter~14]{HrLo15}:
\begin{itemize}
\item $\RR_\infty$ denotes $\RR\cup\{\infty\}$, equipped with the order topology. Moreover, $\RR^+_\infty:=\RR^+\cup\{\infty\}$ denotes the sub-interval of its non-negative elements.
\item We set $\bF:=(F,\RR)$. In particular, $\RR_\infty=\Gamma_\infty(\bF)$ and $\RR^+_\infty=[0,\infty](\bF)$.
\item Let $V$ be an algebraic variety defined over $F$ and $X\subseteq V\times\Gamma_\infty^n$ an $\bF$-definable subset. As a set, let $B_{\bF}(X):=\{p\in S_X(\bF)\mid p\perp^{a}\Gamma\}$. 
\item We  endow $B_{\bF}(V)$ with the topology whose basic open sets are given by finite intersections of sets of the form $\{p\in B_{\bF}(U)\mid g(p)\in\Omega\}$, where $U$ is an open affine subvariety of $V$ defined over $F$, $\Omega\subseteq\RR_\infty$ is an open interval and $g=\val\circ G$ for some $G\in F[U]$.
\item We endow $B_{\bF}(V\times\Gamma_\infty^n)=B_{\bF}(V)\times\RR_\infty^n$ with the product topology, and finally $B_{\bF}(X)\subseteq B_{\bF}(V)\times\RR_\infty^n$ with the subspace topology.
\item Any $\bF$-definable map $f:X\rightarrow Y$ induces a map $B_{\bF}(f):B_{\bF}( X)\rightarrow B_{\bF}(Y)$, which is continuous in case $f$ is the restriction of a regular map between the ambient algebraic varieties.
\item We fix a  maximally complete  algebraically closed  extension $F^{\max}$ of $F$, with $k_{F^{\max}}=(k_F)^{\alg}$ and $\Gamma_{F^{\max}}=\RR$. (Note that such an $F^{\max}$ is uniquely determined up to $\bF$-isomorphism.)
\end{itemize}

The following is an adaptation of Definition~\ref{D:SpecialRetr} to the Berkovich setting.

\begin{definition}
Let $V$ be an 
algebraic variety defined over the valued field $F$, and let $X\subseteq V$ be an $\bF$-definable subset. A continuous map 
$$\widetilde{H}:\RR^+_\infty\times B_{\bF}(X)\rightarrow B_{\bF}(X)$$
is called a \emph{special deformation retraction of $B_{\bF}(X)$ with final image $\bZ$} if the following properties hold:
\begin{itemize}
\item $\widetilde{H}_{\infty}=\mathrm{id}_{B_{\bF}(X)}$
\item $\widetilde{H}_0(B_{\bF}(X))\subseteq\bZ$
\item $\widetilde{H}_t\!\!\upharpoonright_{\bZ}=\mathrm{id}_{\bZ}$ for all $t\in \RR^+_\infty$.
\item For every open subvariety $U$ of $V$ defined over $F$, the set $B_{\bF}(U)\cap B_{\bF}(X)$ is invariant under $\widetilde{H}$.
\item $\widetilde{H}(0,x)=\widetilde{H}(0,\widetilde{H}(t,x))$ for any $x\in B_{\bF}(X)$ and any $t\in \RR^+_\infty$.
\item $\bZ$ is homeomorphic to a piecewise linear subset of $\RR^n$, for some $n\in\NN$.
\item For any $x\in B_{\bF}(X)$ and any $t<\infty$, $\widetilde{H}(t,x)$ is Zariski generic in $X$.
\end{itemize}
\end{definition}

\begin{fact}\label{F:Descent}
Let $V,W$ be algebraic varieties defined over $F$,  and let $X\subseteq V\times\Gamma_\infty^n$ and $Y\subseteq W\times\Gamma_\infty^m$ be $\bF$-definable subsets.
\begin{enumerate}
\item The space $B_{\bF}(V)$ is canonically homeomorphic to the underlying topological space of the Berkovich analytification $V^{\an}$ of $V$.
\item The restriction map $\pi_X:\widehat{X}(F^{\max})\rightarrow B_{\bF}(X)$ is surjective, continuous and closed. If $F=F^{\max}$, it is a homeomorphism.
\item\label{"definable"} Any continuous pro-$\bF$-definable function $h:\widehat{X}\rightarrow\widehat{Y}$ induces a continuous function $\tilde{h}:B_{\bF}(X)\rightarrow B_{\bF}(Y)$ such that $\pi_Y\circ h=\tilde{h}\circ\pi_X$.
\item Let $H:[0,\infty]\times\widehat{X}\rightarrow\widehat{X}$ be an $\bF$-definable special deformation retraction, with final image $H(0,\widehat{X})=Z$. Let $\bZ=\pi_X(Z(F^{\max}))$. Then $H$ induces a special deformation retraction $\widetilde{H}:\RR^+_\infty\times B_{\bF}(X)\rightarrow B_{\bF}(X)$ with final image $\bZ$.\label{Special-preserve}
\end{enumerate}
\end{fact}

\begin{proof}
This is a combination of 14.1.1, 14.1.2, 14.1.3 and 14.1.6 in \cite{HrLo15}, except for the fact in (\ref{Special-preserve}) that $\widetilde{H}$ is a special deformation retraction when $H$ is assumed to 
be an $\bF$-definable special deformation retraction. This is easily checked by hand.
\end{proof}

We will call a  map $\tilde{h}:B_{\bF}(X)\rightarrow B_{\bF}(Y)$ as in part (\ref{"definable"}) of Fact~\ref{F:Descent} \emph{definably induced}. 

\begin{theorem}\label{T:Berk-contract-first}
Let $G$ be an algebraic group defined over $F$, and let $N$ be an $F$-definable stably dominated connected subgroup of $G$. Assume that $N$ is commutative or that $F$ is a valued field of equicharacteristic 0. Then there is a definably induced $N(F)$-equivariant special deformation retraction $\tilde{\rho}: \RR^+_\infty\times B_{\bF}(N)\rightarrow B_{\bF}(N)$ with final image 
$\tilde{\rho}(0,B_{\bF}(N))=\{p_N | \bF\}$. 
\end{theorem}

\begin{proof}
It suffices to apply Fact~\ref{F:Descent}(\ref{Special-preserve}) to Theorem~\ref{T:Contraction-gp} or to Corollary~\ref{C:Contract-char-00}, respectively.
\end{proof}

Similarly, Fact~\ref{F:Descent}(\ref{Special-preserve}) applied to Theorem~\ref{T:Retract-all-char}  or to Corollary~\ref{C:Retract-char-00}, respectively, yields the following result, where we identify (the underlying topological space of) $S^{\an}$ with $B_{\bF}(S)$.

\begin{theorem}\label{T:Berk-contract}
Let $S$ be a semi-abelian variety defined over $F$. Consider the $F$-definable decomposition $0\rightarrow N\rightarrow S\rightarrow\Lambda\rightarrow0$ from Fact~\ref{F:Decomp}. 
Then there is a definably induced  $S(F)$-equivariant special deformation retraction $\tilde{\rho}: \RR^+_\infty\times S^{\an}\rightarrow S^{\an}$ with final image a skeleton $\bf{\Sigma}$. The map $\pi:S\rightarrow\Lambda$ induces a bijection between $\Sigma$ and  $\Lambda(\bF)$, where $\Lambda(\bF)$ is a subset of $\RR^n$ carrying the structure of a piecewise linear abelian group.\qed
\end{theorem}

\begin{remark} Recall that in any NIP theory, in particular in $\acvf$, every $\infty$-definable group $G$ admits a strong connected component $G^{00}$, by a result of Shelah (see, e.g., \cite[Theorem~8.4]{NIPbook}). Assume now that $S$ is a semi-abelian variety defined over a model of $\acvf$,  and let $0\rightarrow N\rightarrow S\rightarrow\Lambda\rightarrow0$ be the decomposition from Fact~\ref{F:Decomp}. As $N=N^{00}$, we have $S/S^{00}\cong\Lambda/\Lambda^{00}$. Thus the failure of strong connectedness of $S$ may be entirely pushed to the piecewise linear world. 

Moreover, if $S=A$ is an abelian variety, in the expansion of $\Gamma$ to a real closed field $\mathcal{R}$, we get $\Lambda\cong S^1(\mathcal{R})^{e}$ for some integer $e$ with $0\leq e\leq \dim(A)$, and $A^{\an}$ is homotopy equivalent to $S^1(\RR)^{e}\cong\Lambda/\Lambda^{00}$. In this sense, the homotopy type of $A^{\an}$ is encoded in the failure of strong connectedness of $A$.
\end{remark}

\section{NIP abelian groups}\label{S:NIP-abgp}
In this section we prove an analogue of our main theorem, but starting with an arbitrary abelian group definable in an NIP theory. The result is of course much weaker, in particular there is no space of generics that could play the role of $\widehat V$. We obtain a directed system of $\infty$-definable subgroups $C(t)$ and instead of those groups having a stably dominated (or generically stable) generic type---which need not exist in general---we ask that they admit a generically stable invariant \emph{measure}. Groups with such a generically stable invariant measure are called \emph{fsg} and can be thought of as the \emph{definably compact} groups in a general NIP theory.

We first recall some definitions, all of which can be found in more details in \cite{NIPbook}. A (Keisler) measure $\mu$ on a definable set $X$ over a model $M$ is a finitely additive probability measure on $M$-definable subsets of $X$. If $X=G$ is a group, then we say that $\mu$ is invariant if $\mu(g\cdot Y)=\mu(Y)$ for every $M$-definable $Y\subseteq X$ and every $g\in G(M)$. A group admitting an invariant measure (over some, or equivalently any, model) is said to be definably amenable. In particular, if $G(M)$ is amenable as a pure group, then it is definably amenable and it follows that any solvable definable group is definably amenable.

Since a type is a special case of a measure (with values in $\{0,1\}$ instead of $[0,1]$), groups with a generically stable invariant type are definably amenable. Many conditions on types translate to measures, and this is the case in particular for generically stability. There are several equivalent definitions of a generically stable measure. We give one which is easy to state. In the following, if $a_1,\ldots,a_n$ are tuples from a model $M$ and $\phi(x)$ is a formula over $M$, we let $Av(\phi(x);a_1,\ldots,a_n)$ denote $|\{i\leq n : M\models \phi(a_i)\}|$.

\begin{definition}
Let $T$ be an NIP theory and $\mu(x)$ a measure over a model $M$. We say that $\mu$ is generically stable if it can be uniformly approximated by finite averages of points in the following sense:

For any formula $\phi(x;y)$ and $\epsilon>0$, there are $a_1,\ldots,a_n\in M^{|x|}$ such that for any $b\in M^{|y|}$, \[| \mu(\phi(x;b)) - Av(\phi(x;b);a_1,\ldots,a_n) | \leq \epsilon.\]
\end{definition}

Note that a generically stable measure $\mu(x)$ is definable in the following sense: given $\phi(x;y)$ and $\epsilon>0$, there are formulas $\theta_i(y,\bar d)$, $i=1,\ldots,n$, that cover $y$-space and such that for any $b,b'$ if for some $i\leq n$, $b,b\models \theta_i(y,\bar d)$, then $|\mu(\phi(x;b))-\mu(\phi(x;b'))|\leq \epsilon$. Indeed, one can take $\bar d$ to consist of the points $a_1,\ldots,a_n$ as in the definition above for $\epsilon/2$ and let the formulas $\theta_i(y,\bar d)$ enumerate the $\phi$-types over $\bar d$. Note in particular that as $\phi(x;y)$ and $\epsilon$ vary, we only need $|T|$-many parameters for $\bar d$, since it is enough to take $\epsilon$ rational.

An $\infty$-definable group $G$ which admits a generically stable invariant measure $\mu(x)$ is said to have fsg (\emph{finitely satisfiable generics}). This condition can be thought of as an abstract version of compactness. For instance, in an o-minimal structure, a definable group has fsg if and only if it is definably compact. Similarly we can state:

\begin{lemma}
	Let $A$ be an abelian variety definable in a model of $\acvf$, then $A$ has fsg.
\end{lemma}
\begin{proof} Let $0\rightarrow N\rightarrow A\rightarrow\Lambda\rightarrow0$ be the decomposition from Fact~\ref{F:Decomp}. The group $N$ is stably dominated (hence generically stable), so it has fsg. By \cite[Corollary~6.19]{HrRi19}, $\Lambda$ is a definably compact, $\Gamma$-internal group, hence $\Lambda$ has fsg by \cite[Theorem~8.1]{NIP1}. 
By Proposition~4.5 of that same paper, an extension of a group with fsg by a group with fsg has fsg. Hence $A$ has fsg.
\end{proof}


If $\mu(x)$ and $\eta(x)$ are two generically stable (Keisler) measures on a definable group $G$, we can define the convolution $\mu \ast \eta (x)$ by \[\mu\ast \eta(\phi(x))=\int_y \mu(\phi(x\cdot y))d\lambda(y).\] We refer to \cite[Section~7.4]{NIPbook} for explanation of why this integral makes sense. Then $\mu \ast \eta$ is again a generically stable measure. Note that if $G$ is abelian, then for any $g\in G$, 
\begin{itemize}
\item[(A)] 
$g\cdot (\mu\ast\eta) = (g\cdot \mu)\ast \eta = \mu \ast (g\cdot \eta).$
\end{itemize}

\begin{proposition}\label{P:stab-Keisler}
Let $T$ be NIP and let $G$ be a definable abelian group. Then there is a pro-definable set $S$ and $\infty$-definable subgroups $N(t)$, for $t\in S$ forming a directed system (any small family has an upper bound) and such that $\bigcup_{t\in S} N(t) = G$ and each $N(t)$ stabilizes a generically stable measure on $G$.
\end{proposition}

\begin{proof}
As $G$ is abelian, it is definably amenable. Let $M$ be a $|T|^+$-saturated model and let $\mu_M$ be a $G$-invariant measure over $M$. We can extend $\mu_M$ to a global measure $\mu$ which is $G(M)$-invariant and generically stable (see \cite[Lemma~7.6]{NIP2} or \cite[Proposition~3.4]{Finding}). Let $N_\mu= \{g\in G(\monster) \mid   g\cdot \mu = \mu\}$ be the stabilizer of $\mu$. Then $N_\mu$ is a subgroup of $G$ containing $G(M)$. By definability of $\mu$, $N_\mu$ is $\infty$-definable over a set of size $|T|$. Let $t_0$ be an enumeration of the parameters needed to define $\mu$ and write $\mu=\mu_{t_0}$. Then also $N_\mu$ is defined over $t_0$ and we can write $N_\mu=N(t_0)$. Let $S=\tp(t_0)$. Then $S$ is pro-definable and for every tuple $t\in S$, $\mu_t$ is a well-defined generically stable measure (where $\mu_t$ is defined over $t$ using the same definition scheme as $\mu$ over $t_0$). Then also $N_{\mu_t}=N(t)$ is defined over $t$ the same way $N_\mu$ is defined over $t_0$. We have constructed a pro-definable family of $\infty$-definable subgroups of $G$.

Since $N(t_0)$ contains $G(M)$, by compactness, for any (small) model $M'$, there is $t'\equiv t_0$ such that $N(t')$ contains $G(M')$. This proves that $\bigcup_{t\in S} N(t)=G$.

It remains to show that the family is directed. Let $(r_i)_{i\in I}$ be a small family of elements of $S$ and write $\mu_i=\mu_{r_i}$. For any finite $I_0\subseteq I$, let $\mu_{I_0}= \ast_{i\in I_0} \mu_i$. As $G$ is abelian, this product is independent of the order of the factors. By (A), each $N(r_i)$, $i\in I_0$ is in $Stab(\mu_{I_0})$. Fix a model $M'$ over which all the $\mu_i$'s are defined. Let $\mu$ be a limit of the $\mu_{I_0}$ over $M'$ along an ultrafilter containing $\{I_1:I_1\supseteq I_0\}$ for each finite $I_0\subseteq I$. Then $\mu$ is $M'$-invariant and its stabilizer contains each $N(r_i)$, $i\in I$. As above, we can extend $\mu$ to a global measure $\lambda$ which is generically stable and whose stabilizer still contains all the $N(r_i)$'s. Now $\lambda$ is defined over some model of size $|T|$ and hence there is an automorphism $\sigma$ such that $\sigma(\lambda)=: \lambda_0$ is generically stable over $M$. Let $t_0$ be as in the first paragraph of the proof. Then $N(t_0)$ contains $G(M)$ and as $\lambda_0$ is finitely satisfiable in $M$, we have $\lambda_0(N(t_0))=1$. This implies that $Stab(\lambda_0)\leq N(t_0)$. So $Stab(\lambda)\leq N(\sigma^{-1}(t_0))$, with $\sigma^{-1}(t_0)\in S$. Therefore each $N(t_i)$ is a subgroup of $N(\sigma^{-1}(t_0))$, which finishes the proof.
\end{proof}

\begin{corollary}
If $G$ is an abelian group with no indiscernible linearly ordered family of $\infty$-definable groups, then $G$ has fsg.
\end{corollary}
\begin{proof}
By the previous proposition and since the family $C(t)$ is directed, it must be that $G=C(t)$ for some $t$. Then $G$ stabilizes a generically stable measure and hence $G$ has fsg.
\end{proof}

\bigskip
\bibliographystyle{plain}
\bibliography{Retract}
\end{document}